\documentclass[11pt]{amsart}
\usepackage{graphicx}
\usepackage{color}
\usepackage{amssymb}
\usepackage{amsmath}
\usepackage{epstopdf}
\usepackage{epsfig}
\usepackage{psfrag}
\numberwithin{figure}{section}
\newtheorem{theorem}{Theorem}[section]
\newtheorem{corollary}[theorem]{Corollary}
\newtheorem{lemma}[theorem]{Lemma}
\newtheorem{proposition}[theorem]{Proposition}
\newtheorem{definition}[theorem]{Definition}
\newtheorem{remark}[theorem]{Remark}
\newtheorem{prop}[theorem]{Proposition}

\title{Stabilizations of Heegaard splittings of graph manifolds}
\author{Ryan Derby-Talbot}
\date{\today}

\begin{document}
\maketitle

\begin{abstract}
We show that after one stabilization, a strongly irreducible Heegaard splitting of suitably large genus of a graph manifold is isotopic to an amalgamation along a modified version of the system of canonical tori in the JSJ decomposition. As a corollary, two strongly irreducible Heegaard splittings of a graph manifold of suitably large genus are isotopic after at most one stabilization of the higher genus splitting.
\end{abstract}

\newcommand{\hs}{$V \cup_S W$}
\newcommand{\hsone}{$V_1 \cup_{S_1} W_1$}
\newcommand{\hstwo}{$V_2 \cup_{S_2} W_2$}
\newcommand{\hhs}{$V' \cup_{S'} W'$}
\newcommand{\pq}{$P \cup_{\Sigma} Q$}
\newcommand{\ppq}{$P' \cup_{\Sigma'} Q'$}
\psfrag{S}{$S$}
\psfrag{F}{$F$}
\psfrag{L}{$L$}
\psfrag{Tp}{$T_1'$}
\psfrag{Tq}{$T_2'$}
\psfrag{T1}{$T_1$}
\psfrag{T2}{$T_2$}
\psfrag{Te}{$T_{e_1}$}
\psfrag{Tepr}{$T_{e_1}'$}
\psfrag{X1}{$X_1$}
\psfrag{X2}{$X_2$}
\psfrag{e1}{$e_1$}
\psfrag{e2}{$e_2$}
\psfrag{e3}{$e_3$}
\psfrag{g1}{$\gamma_1$}
\psfrag{g2}{$\gamma_2$}
\psfrag{w1}{$\omega_2$}
\psfrag{w2}{$\omega_3$}
\psfrag{l1}{$\lambda_1$}
\psfrag{l2}{$\lambda_2$}
\psfrag{B}{$B$}
\psfrag{d1}{$\delta_1$}
\psfrag{d2}{$\delta_2$}
\psfrag{d3}{$\delta_3$}
\psfrag{D1}{$D_1$}
\psfrag{t1}{$\tau_1$}
\psfrag{t2}{$\tau_2$}
\psfrag{G}{$G$}
\psfrag{B1}{$A_1$}
\psfrag{B2}{$A_2$}
\psfrag{F0}{$F_0$}
\psfrag{VR}{$V_i$}
\psfrag{Vj}{$V_j$}
\psfrag{WR}{$W_i$}
\psfrag{Wj}{$W_j$}
\psfrag{p1}{{\small $\partial_1 X_j$}}
\psfrag{p2}{{\small $\partial_2 X_i$}}
\psfrag{pp2}{{\tiny $\partial_2 X_j - F_0$}}
\psfrag{pp1}{{\tiny $\partial_1 X_i - F_0$}}
\psfrag{Xi}{$X_i$}
\psfrag{(X)}{$X_i$}
\psfrag{Xj}{$X_j$}
\psfrag{Fm}{$F_m$}
\psfrag{V}{$V$}
\psfrag{W}{$W$}
\psfrag{alpha}{$\alpha$}
\psfrag{Vp}{$V'$}
\psfrag{Wp}{$W'$}
\psfrag{N}{$N$}
\psfrag{pN}{$\partial N$}
\psfrag{BR}{{\tiny boundary compressing disk for $S$}}
\psfrag{st}{stabilization}
\psfrag{S0}{$Z_i$}
\psfrag{A}{$A$}
\psfrag{Q}{$Q$}
\psfrag{t}{$\tau$}
\psfrag{tp}{$\tau'$}
\psfrag{V2}{$V' \cap X_1$}
\psfrag{di}{disk}
\psfrag{du}{dual tube}
\psfrag{pQ}{$\partial Q$}
\psfrag{R}{$\partial X$}
\psfrag{pX}{$\partial X$}
\psfrag{Xprime}{$X'$}
\psfrag{X}{$X$}
\psfrag{dui}{dual (if necessary)}
\psfrag{D2}{$D_2$}
\psfrag{un}{{\small $\begin{array}{l} \mbox{ components of} \\ W' \cap (\partial Y - \partial M) \end{array}$}}
\psfrag{u1}{{\small unblocked}}
\psfrag{pS}{$S'$}
\psfrag{sp}{{\small spanning disks}}
\psfrag{ER}{{\small exceptional fiber}}
\psfrag{c}{$\cong$}
\psfrag{bb}{{\small blocked annulus}}
\psfrag{bl}{{\small components of $V' \cap (\partial Y - \partial M)$}}
\psfrag{i}{(i)}
\psfrag{ii}{(ii)}
\psfrag{iii}{(iii)}
\psfrag{Th}{$\Theta$}
\psfrag{TR}{$\Theta$}
\psfrag{tr}{{\small trivial tube}}
\psfrag{U}{$U$}
\psfrag{H}{$H$}
\psfrag{Y1}{$Y_1$}
\psfrag{a}{{\small annuli}}
\psfrag{Tt}{$T^2 \times I$}
\psfrag{Xm}{$X''$}
\psfrag{Up}{$U'$}
\psfrag{AR}{$A'$}
\psfrag{or}{or}
\psfrag{T}{$T$}
\psfrag{Tm}{$T'$}
\psfrag{MT}{$T'$}
\psfrag{UI}{$U_i$}
\psfrag{Oi}{$\Upsilon_i$}
\psfrag{Oj}{$\Upsilon_j$}
\psfrag{O}{$\Upsilon$}
\psfrag{D}{$D$}
\psfrag{E}{$E$}
\psfrag{CR}{{\small $\begin{array}{l} \mbox{ components of} \\ \ \ \ W' \cap X \end{array}$}}
\psfrag{A1}{$A_1$}
\psfrag{A2}{$A_2$}
\psfrag{D0}{$D_0$}

\section{Introduction}

Reidemeister and Singer proved in the 1930s that any two Heegaard splittings of the same 3-manifold become isotopic after stabilizing each splitting some number of times (\cite{Reidemeister}, \cite{Singer}). It has been a long-standing problem to determine the precise number of stabilizations needed to do so. This problem is made challenging by the fact that stabilized Heegaard splittings lose many geometric properties that non-stabilized Heegaard splittings possess, making their analysis difficult. Many attempts to understand the number of necessary stabilizations have exploited the underlying structure of the 3-manifold (see {\em e.g.}~\cite{RS}, \cite{Schultens1} and \cite{Waldhausen}). This paper is no exception. Here, we investigate the question of how many stabilizations are needed to make two Heegaard splittings of a graph manifold isotopic. Our main results are the following (see below for defintions). 

\begin{theorem}
\label{maintheorem}
Let $M$ be a closed totally orientable graph manifold with $\Theta$ a system of canonical tori in the JSJ decomposition, and let $\Theta'$ be an amalgamatable modification of $\Theta$. Suppose \hs \ and \pq \ are Heegaard splittings of $M$ such that \hs \ either is an amalgamation along $\Theta'$ or is strongly irreducible, and \pq \ is an amalgamation along $\Theta'$. Then \hs \ and \pq \ are isotopic after at most one stabilization of the larger genus splitting.
\end{theorem}

\begin{corollary}
\label{maincorollary}
If $M$ has two strongly irreducible Heegaard splittings, at least one of which has genus greater than or equal to $a(M)$, then the two splittings are isotopic after at most one stabilization of the larger genus splitting. 

\end{corollary}

An amalgamatable modification of $\Theta$ (see Definition~\ref{amalgamatablemodificationdefinition}) is a disjoint union of incompressible tori in $M$ containing $\Theta$ that has the appropriate separating properties to yield an amalgamation. The amalgamation genus of $M$, $a(M)$ (see Definition~\ref{amalgamationgenusdefinition2}), is the minimum possible genus of a Heegaard splitting of $M$ which is an amalgamation along an amalgamatable modification of $\Theta$. 

We in fact prove more general versions of Theorem~\ref{maintheorem} and Corollary~\ref{maincorollary} including the case that $M$ has boundary. The statements of these results are slightly more technical and are stated as Theorem~\ref{maintheorem1} and Corollary~\ref{maincorollary1}, respectively. The main difference is that we need to consider $M$ as a triple $(M, \partial_1 M, \partial_2 M)$, and take the amalgamation genus as a quantity $a(M, \partial_1 M, \partial_2 M)$ depending on a partition of the boundary components. 

It is not difficult to see that two Heegaard splittings of $M$ which are amalgamations along the same amalgamatable modification of $\Theta$ are isotopic after at most one stabilization of the higher genus splitting (see Corollary~\ref{graphmanifoldsamalgamationlemma}). It is the case that at least one of the splittings is strongly irreducible that is the main focus of this paper. 

\begin{remark}
\label{firstremark}
\rm Theorem~\ref{maintheorem} and Corollary~\ref{maincorollary} along with Theorem 1.3 in \cite{Schultens2} imply that two irreducible Heegaard splittings of a closed, totally orientable graph manifold $M$ are isotopic after at most one stabilization of the larger genus splitting, unless:

\begin{itemize}
\item[1.] The two Heegaard splittings are both strongly irreducible with genus less than $a(M)$, or

\item[2.] One of the splittings is an amalgamation along some incompressible surface and is not isotopic to an amalgamation along $\Theta'$, where $\Theta'$ is an amalgamatable modification of the system of canonical tori in the JSJ decomposition of $M$.
\end{itemize}

Indeed, there are graph manifolds where both situations are possible: see Section~\ref{sec:section11} in this paper for examples of the first situation and \cite{Schultens2} for examples of the second. An analogous remark can be made in the case $M$ has boundary, considering a partition $(M, \partial_1 M, \partial_2 M)$. 

\end{remark}

While the arguments are somewhat delicate, the strategy of the proof of Theorem~\ref{maintheorem} is straightforward. In Section~\ref{sec:background} we provide the basic definitions needed for Heegaard splittings, Seifert fibered spaces and graph manifolds. In Section~\ref{amalgamationsandamalgamationgenussection}, we establish properties of Heegaard splittings which are amalgamations along incompressible surfaces, paying particular attention to the case that the manifold in question is a graph manifold and the incompressible surface used in amalagamation contains the system of canonical tori in the JSJ decomposition of the manifold. In Section~\ref{sec:section4}, we determine the precise structure of a strongly irreducible Heegaard splitting of a graph manifold, drawing mostly from results in \cite{Schultens2}. There it was shown that a strongly irreducible Heegaard splitting of a graph manifold can be taken to be horizontal or vertical in all but one of the Seifert fibered components of the manifold. In Section~\ref{sec:section5}, we prove a special case of Theorem~\ref{maintheorem} in which there are no Seifert fibered components where the Heegaard surface is horizontal. In Section~\ref{sec:horizontalcomponents} we deal with the remaining situation, showing that the Seifert fibered components where the Heegaard surface is horizontal use one stabilization to produce many more, which then can be dispersed throughout the manifold. In Section~\ref{sec:section7} we discuss the precise way the stabilizations are isotoped in the manifold. In Sections~\ref{sec:section8} and \ref{sec:section9} we argue that the stabilizations produced by the horizontal pieces can be used in manifolds homeomorphic to $T^2 \times I$ and in vertical pieces. In Section~\ref{sec:proofofthemainresults}, we prove Theorem~\ref{maintheorem} by showing that if there are enough stabilizations, a quantity that depends only on the genus of the splitting, then the stabilized Heegaard splitting becomes an amalgamation along an amalgamatable modification of the system of canonical tori in the JSJ decomposition. In Section~\ref{sec:section11} we provide examples of both large and small genus Heegaard splittings of graph manifolds.

The author wishes to thank both Cameron Gordon and Tsuyoshi Kobayashi for many helpful discussions.

\section{Background}
\label{sec:background}

For basic terminology involving 3-manifolds, see \cite{Hempel}, \cite{Jaco} or \cite{Rolfsen}.

\subsection{Heegaard splittings}

\begin{definition} 
\rm Let $F$ be a closed orientable surface. A {\em compression body} $V$ is obtained from $F \times I$ by adding 2-handles to $F \times \{0\}$ and capping off any resulting 2-sphere components with 3-balls. The surface $F \times \{1\}$ is denoted $\partial_+ V$, and we set $\partial_- V = \partial V - \partial_+ V$. If $\partial_- V = \emptyset$, then $V$ is called a {\em handlebody}.
\end{definition}

\begin{definition}
\rm Let $M$ be a compact orientable 3-manifold. A {\em partition} of $M$ is a triple $(M, \partial_1 M, \partial_2 M)$ such that $\partial_1 M$ and $\partial_2 M$ consist of components of $\partial M$, $\partial_1 M \cap \partial_2 M = \emptyset$, and $\partial M = \partial_1 M \sqcup \partial_2 M$. We call the triple $(M, \partial_1 M, \partial_2 M)$ a {\em partitioned} 3-manifold.
\end{definition}

\begin{definition}
\rm A {\em Heegaard splitting} \hs \ of $(M, \partial_1 M, \partial_2 M)$ is a decomposition of $M$ into 2 compression bodies $V$ and $W$ such that $\partial_- V = \partial_1 M$ and $\partial_- W = \partial_2 M$, and such that $M$ is obtained from the identification of $\partial_+ V$ and $\partial_+W$ via some homeomorphism. The surface $S = \partial_+ V = \partial_+ W$ is called the {\em splitting surface} or the {\em Heegaard surface}.
\end{definition}

If $M$ is closed, the partition $(M, \emptyset, \emptyset)$ is written as $M$, so that \hs \ is simply referred to as a Heegaard splitting of $M$. The {\em genus} of a Heegaard splitting \hs \ of $(M, \partial_1 M, \partial_2 M)$ is the genus of $S$. It is known that every 3-manifold has a Heegaard splitting (\cite{Moise}). Also, note that a Heegaard splitting of $(M, \partial_2 M, \partial_1 M)$ is a Heegaard splitting of $(M, \partial_1 M, \partial_2 M)$ with the labels of the compression bodies $V$ and $W$ switched. 

Suppose \hs \ is Heegaard splitting of $(M, \partial_1 M, \partial_2 M)$ and $\alpha$ is an arc properly embedded in $V$ and parallel into $\partial_+ V$. Let $V'$ and $W'$ be obtained from $V$ and $W$ by attaching a 1-handle to $W$ such that $\alpha$ is its core, and removing that 1-handle from $V$. The fact that $\alpha$ is boundary parallel implies that both $V'$ and $W'$ are compression bodies. The resulting Heegaard splitting \hhs \ of $(M, \partial_1 M, \partial_2 M)$ is called a {\em stabilization} of \hs . Note that the genus of \hhs \ is one greater than the genus of \hs . 
 
\begin{figure}[h]
\centering
\includegraphics[width=4in]{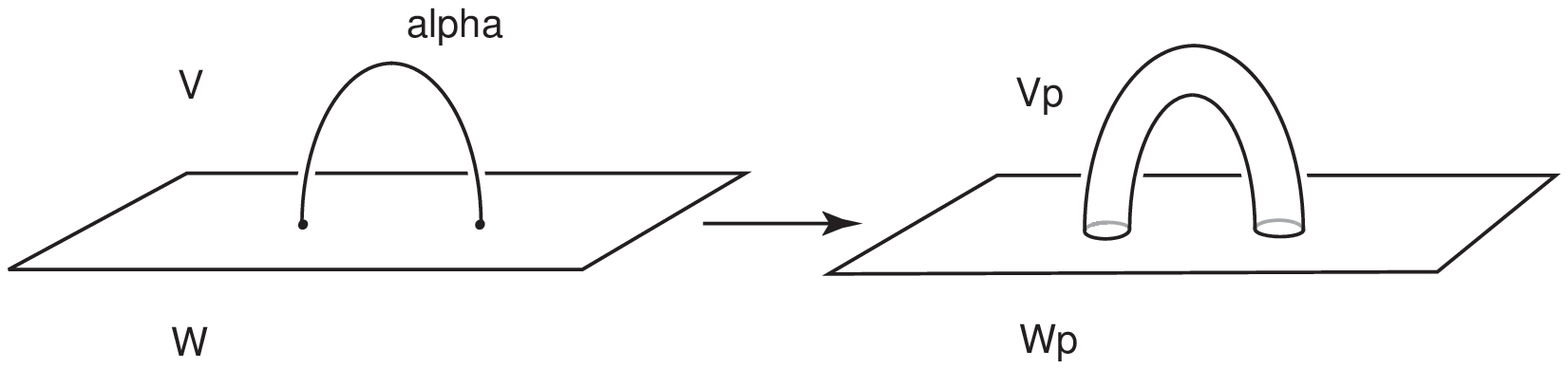}
\caption{A stabilization of \hs .} 
\label{fig:stabilization}
\end{figure}

It is a nice exercise to show that a Heegaard splitting \hs \ is stabilized ({\em i.e.}~a stabilization of another Heegaard splitting) if and only if there exist essential disks $D \subset V$ and $E \subset W$ such that $|\partial D \cap \partial E| = 1$. 

\begin{definition}
\rm A Heegaard splitting \hs \ of $(M, \partial_1 M, \partial_2 M)$ is {\em reducible} if there exist essential disks $D \subset V$ and $E \subset W$ such that $\partial D = \partial E$. Otherwise, \hs \ is {\em irreducible}. 

The splitting \hs \ is {\em weakly reducible} if there exist essential disks $D \subset V$ and $E \subset W$ such that $\partial D \cap \partial E = \emptyset$. Otherwise, \hs \ is strongly irreducible.
\end{definition}

A reducible Heegaard splitting is clearly weakly reducible. A result of Haken \cite{Haken} shows that a Heegaard splitting of a reducible manifold is reducible. In light of this fact we assume henceforth that $M$ is irreducible. A strongly irreducible Heegaard splitting has the useful property that it can be isotoped to intersect an incompressible surface in curves which are essential in both surfaces (see {\em e.g.}~\cite{Schultens3}). 

\begin{remark}
\label{startingpositionremark}
\rm Suppose that \hs \ is a strongly irreducible Heegaard splitting of $(M, \partial_1 M, \partial_2 M)$ and that $F$ is a disjoint union of incompressible tori. After performing the aforementioned isotopy and then isotoping away boundary parallel components, $F$ intersects both $V$ and $W$ in a disjoint union of essential annuli (an {\em essential annulus} is defined in Section~\ref{amalgamationsandamalgamationgenussection}). This can be done so that $|S \cap F|$ is minimal. We shall assume henceforth that any strongly irreducible Heegaard splitting of $(M, \partial_1 M, \partial_2 M)$ is isotoped in such a fashion. 
\end{remark}

\subsection{Seifert fibered spaces} 

Let $B$ be a surface and set $X' = B \times S^1$. Let $e_1, \ldots, e_k$ be a (possibly empty) disjoint union of points on $B$ (called {\em exceptional points}), and let $U_i$ be the solid torus in $X'$ formed by taking $N(e_i) \times S^1$, where $N(e_i)$ denotes a closed neighborhood of $e_i$, disjoint from $N(e_j)$ for $j \neq i$. Let $\mu_i\subset \partial U_i$ bound a meridian disk of $U_i$. Consider $B = B \times \{point\} \subset X'$. Let $X_0 = \overline{X' - \cup_{i=1}^k U_i}$ and let $B_0 = X_0 \cap B$. For each boundary component in $X_0$ corresponding to $\partial U_i$, let $c_i = B_0 \cap \partial U_i$ and let $t_i = \{point\} \times S^1 \subset \partial U_i$. Form the 3-manifold $X$ by re-attaching $U_i$ to $X_0$ by gluing $\mu_i$ to a curve corresponding to $c_i^{\alpha_i} t_i^{\beta_i}$. Once the meridian is attached, there is a unique way to attach the rest of $U_i$ to $X_0$ up to homeomorphism. We say that $X$ is a {\em Seifert fibered space} with {\em base} $B$, and {\em exceptional fibers} $f_1, \ldots, f_k$ where $f_i$ is the core of $U_i$ (assuming $|\alpha_i| > 1$). Every other curve $\{ point \} \times S^1$ in $X$ is a {\em regular fiber}. If both $B$ and $X$ are orientable, we say $X$ is {\em totally orientable}. For the purposes of this paper, we shall always assume that $X$ (and hence $B$) has boundary. Note that each component of $\partial X$ is a torus.

All of the data describing $X$ as a Seifert fibered space with boundary up to homeomorphism can be written in the notation 

\begin{displaymath}
\left(g; m; \frac{\beta_1}{\alpha_1}, \ldots, \frac{\beta_k}{\alpha_k} \right)
\end{displaymath}

\noindent where $g$ is the genus of $B$, $m$ is the number of boundary components of $B$, and $\alpha_i$, $\beta_i$ are the coordinates in the $\{c_i, t_i\}$ basis of the attaching curve for the meridian of $U_i$. The number $\alpha_i$ is called the {\em multiplicity} of the exceptional fiber $f_i$. The fraction $\frac{\beta_i}{\alpha_i}$ is unique up to adding integers.

\subsection{Graph manifolds}

\begin{definition}
\label{graphmanifolddefinition}
\rm A {\em graph manifold} $M$ is a 3-manifold obtained by gluing a disjoint union of Seifert fibered spaces together along boundary tori. 
\end{definition} 

Let $\Theta$ denote the system of canonical tori in the JSJ decomposition of $M$. For a graph manifold, each component of $M$ cut along $\Theta$ is a Seifert fibered space, referred to as a {\em Seifert fibered component} of $M$. We say that $M$ is {\em totally orientable} if it is orientable and each Seifert fibered component of $M$ is totally orientable. 

\begin{remark}
\label{graphmanifoldsremark}
\rm It is sometimes useful to consider the manifolds $N(\Theta)$ in a graph manifold $M$, where $N(\Theta)$ denotes a closed regular neighborhood of $\Theta$ in $M$. Assuming the neighborhood $N(\Theta)$ meets every Seifert fibered component of $M$ cut along $\Theta$ in regular fibers, then $\overline{M - N(\Theta)}$ is a disjoint union of Seifert fibered spaces homeomorphic to the components of $M$ cut along $\Theta$. These components are still referred to as the Seifert fibered components of $M$. Each component of $N(\Theta)$ is homeomorphic to $T^2 \times I$, and is called an {\em edge manifold} of $M$. 
\end{remark}

We will only use this description of graph manifolds in Section~\ref{sec:section4}. Otherwise, when stating $M$ is a graph manifold we mean in the context of Definition~\ref{graphmanifolddefinition}.

\section{Amalgamations and amalgamation genus}
\label{amalgamationsandamalgamationgenussection}

\subsection{Amalgamations}

Let $A$ be an annulus properly embedded in a compression body $V$. We say $A$ is {\em essential} in $V$ if $A$ is incompressible in $V$ and not boundary parallel. A {\em spanning arc} of $A$ is a properly embedded arc $a$ in $A$ such that $\partial a$ lies in different components of $\partial A$. Hence $A$ cut along a spanning arc is a disk.

Let $A$ be an essential annulus in a compression body such that $\partial A \subset \partial_+ V$. A {\em spanning disk} for $A$ in $V$ is a disk $D$ in $V$ cut along $A$ such that $\partial D = a \cup b$,  where $a$ is a spanning arc of $A$ and $b \subset \partial_+V$. It follows that $D$ is an essential disk in $V$ cut along $A$.  

\begin{proposition}
\label{spanningdisk}
An essential annulus $A$ in a compression body $V$ with $\partial A \subset \partial_+ V$ has a spanning disk. 
\end{proposition}

\begin{proof}
Let $\Delta$ be a complete system of meridian disks for $V$. Then $V$ cut along $\Delta$ is either a 3-ball if $V$ is a handlebody, or a product $\partial_- V \times I$ otherwise. By a standard innermost disk, outermost arc argument, we can eliminate inessential curves and arcs of intersection in $A \cap \Delta$ from $A$. If $V$ cut along $\Delta$ is a 3-ball, then $A \cap \Delta \neq \emptyset$ since there are no incompressible annuli in a 3-ball. An outermost disk component of a component of $\Delta$ cut along $A$ is thus a spanning disk for $A$. 

Assume therefore that $V$ cut along $\Delta$ is a product $\partial_- V \times I$ so that $\partial A \subset \partial_- V \times \{ 1 \}$. Let $C$ be a union of simple closed curves intersecting transversely on $\partial_- V$ such that $\partial_- V$ cut along $C$ is a disjoint union of disks. Again by an innermost disk, outermost arc argument we may assume that $A$ intersects $C \times I$ in spanning arcs of $A$. Then $\partial_- V \times I$ cut along $C \times I$ is a union of 3-balls, so as before, it must be the case that $A$ meets some annulus $c \times I$ nontrivially, where $c$ is a simple closed curve in $C$. As $\partial A \subset \partial_- V \times \{ 1 \}$, $A \cap (c \times I)$ cuts off disk components from $c \times I$. An outermost such disk is thus a spanning disk of $A$. 
\end{proof}

\begin{definition}
\label{mutuallyseparatingdefinition}
\rm An embedded closed surface $F$ in a 3-manifold $M$ is called {\em mutually separating} if $M$ cut along $F$ is two (possibly disconnected) 3-manifolds $M_1$ and $M_2$ such that a regular neighborhood of each component of $F$ meets both $M_1$ and $M_2$ nontrivially.
\end{definition}

\begin{remark}	
\rm The condition that $F$ is mutually separating in $M$ is equivalent to the statement that $F$ represents the trivial element in $H_2(M, \partial M ; \mathbb{Z}/2 \mathbb{Z})$.
\end{remark}

Let $F$ be a mutually separating incompressible surface in $(M, \partial_1 M, \partial_2 M)$, and let $X_1, \ldots, X_n$ be the components of $M$ cut along $F$. The fact that $F$ is mutually separating allows us to choose partitions $(X_i, \partial_1 X_i, \partial_2 X_i)$ for each $1 \leq i \leq n$ so that if $X_i \cap X_j \neq \emptyset$ for any $i \neq j$, then either $F \cap X_i \subset \partial_1 X_i$ and $F \cap X_j \subset \partial_2 X_j$, or $F \cap X_i \subset \partial_2 X_i$ and $F \cap X_j \subset \partial_1 X_j$. 
 
Let $V_i \cup_{S_i} W_i$ be a Heegaard splitting of $(X_i, \partial_1 X_i, \partial_2 X_i)$ for each $1 \leq i \leq n$. Suppose that $F_0 =  \partial_- V_i \cap \partial_- W_j = \partial_1 X_i \cap \partial_2 X_j$ for some $i \neq j$. Then $V_i$ can be thought of as $F_0 \times I \sqcup (\partial_1 X_i - F_0) \times I$ with 1-handles attached along $F_0 \times \{ 1 \} \sqcup (\partial_1 X_i - F_0) \times \{ 1 \}$. The compression body $W_j$ is defined similarly. Map $F_0 \times \{t\}$ to $F_0 \times \{0\}$ in both $V_i$ and $W_j$ such that the disks of attachment of the 1-handles in $V_i$ and in $W_j$ map disjointly on $F_0$. This yields compression bodies $V = V_j \cup (\partial_1 X_i - F_0) \times I \cup \{\mbox{1-handles in $V_i$}\}$, and $W = W_i \cup (\partial_2 X_j - F_0) \times I \cup \{\mbox{1-handles in $W_j$}\}$. Note that $\partial_- V = \partial_1 X_j \sqcup (\partial_1 X_i - F_0)$, and $\partial_- W = \partial_2 X_i \sqcup (\partial_2 X_j - F_0)$, thus \hs \ is a Heegaard splitting of $(X_i \cup X_j, \partial_1 X_j \sqcup (\partial_1 X_i - F_0), \partial_2 X_i \sqcup (\partial_2 X_j - F_0))$. The above procedure can be applied simultaneously to all components of $M$ cut along $F$.

\begin{definition}
\rm A Heegaard splitting of $(M, \partial_1 M, \partial_2 M)$ obtained in the above manner is called an {\em amalgamation along $F$}.
\end{definition}
\begin{figure}[h] 
   \centering
   \includegraphics[width=5.5in]{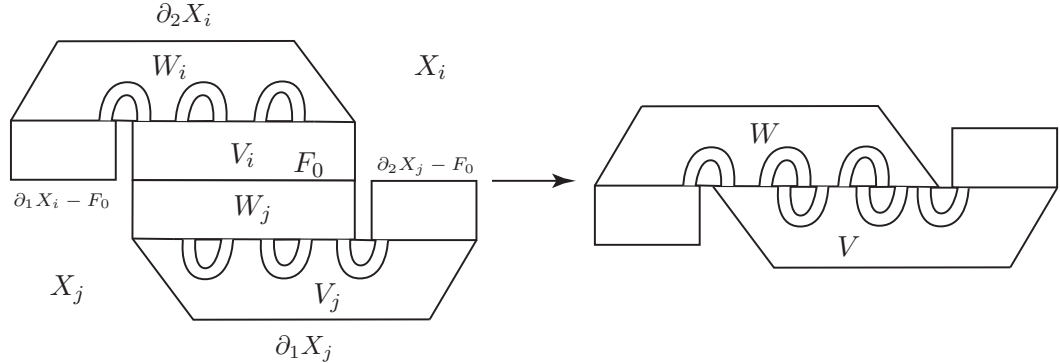}
   \caption{Forming an amalgamation along $F_0$.}
   \label{fig:amalgamationtraditional}
\end{figure}

\begin{remark}
\rm The reverse of the above process is called {\em untelescoping} the Heegaard splitting \hs , and is a useful procedure when dealing with amalgamations along incompressible surfaces (see \cite{ST}).
\end{remark}

\begin{remark}
\rm Note that if $F$ is an incompressible surface which is mutually separating in $M$, a Heegaard splitting \hs \ is an amalgamation along $F$ if $S$ can be simultaneously compressed in $V$ and $W$ to become isotopic to $F \cup L$, where $L$ is a (possibly empty) subset of components of $\partial M$.
Specifically, if the components $X_1, \ldots, X_n$ of $M$ cut along $F$ are partitioned as above so that $F \cap X_i \subset \partial_{\varepsilon_i} X_i$, $\varepsilon_i = 1$ or $2$, then $L = \cup_{i=1}^n (\partial_{\varepsilon_i} X_i - F)$. 
\end{remark}

\begin{figure}[h] 
   \centering
   \includegraphics[width=5in]{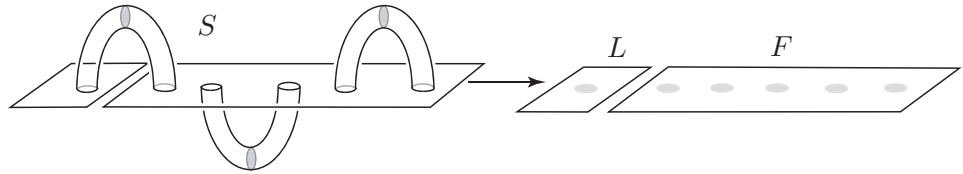} 
   \caption{An amalgamation along $F$.}
   \label{fig:amalgamation}
\end{figure}

It is known that an irreducible Heegaard splitting of $(M, \partial_1 M, \partial_2 M)$ is either strongly irreducible or an amalgamation along some incompressible surface (see \cite{CG} for the closed case and {\em e.g.}~\cite{Kobayashi} for the case that $M$ has boundary).

\subsection{Amalgamation genus}
\label{amalgamationgenussubsection}

\begin{definition}
\rm The {\em Heegaard genus} of $(M, \partial_1 M, \partial_2 M)$, denoted $g_h(M, \partial_1 M, \partial_2 M)$,  is the minimum genus of all Heegaard splittings of $(M, \partial_1 M, \partial_2 M)$. 
\end{definition}

\begin{definition}
\label{amalgamationgenusdefinition1}
\rm The {\em amalgamation genus} of $(M, \partial_1 M, \partial_2 M)$ with respect to a mutually separating incompressible surface $F$, denoted $a(M, \partial_1 M, \partial_2 M, F)$, is the minimum genus of all Heegaard splittings of $(M, \partial_1 M, \partial_2 M)$ which are amalgamations along $F$. 
\end{definition}

\begin{remark}
\label{twopartitionsremark}
\rm Given $(M, \partial_1 M, \partial_2 M)$ containing a mutually separating incompressible surface $F$, the components $X_1, \ldots, X_n$ of $M$ cut along $F$ have 2 natural choices of partitions that can be used to construct an amalgamation along $F$. For some $X_i$, choose either $\partial_1 X_i = (\partial_1 M_1 \cap X_i) \sqcup (F \cap X_i)$ and $\partial_2 X_i = \partial_2 M \cap X_i$, or $\partial_1 X_i = \partial_1 M \cap X_i$ and $\partial_2 X_i = (\partial_2 M \cap X_i) \sqcup (F \cap X_i)$. Having chosen $(X_i, \partial_1 X_i, \partial_2 X_i)$ for this $i$ then determines $(X_j, \partial_1 X_j, \partial_2 X_j)$ for $j \neq i$, via the definition of amalgamation along $F$. As there are two choices for a partition of $X_i$, this gives two choices of sets of partitions for all the $X_i$, $1 \leq i \leq n$. This choice has bearing on the genus of an amalgamation of $(M, \partial_1 M, \partial_2 M)$ along $F$. 

For example, let $X_1$ and $X_2$ be 3-manifolds homeomorphic to $T^2 \times I$ and label $\partial X_i = T_i \cup T_i'$, $i = 1, 2$. Obtain $M$ by identifying $T_1 = T_2 = T$ via any homeomorphism, and take the partition $(M, T_1', T_2')$. There are two natural choices for amalgamations of $(M, T_1', T_2')$ along $T$, depending on the choice of partitions of $X_1$ and $X_2$. If we take the partitions $(X_1, T_1, T_1')$ and $(X_2, T_2', T_2)$, then the resulting amalgamation along $T$ has genus $1$. If, however, we choose partitions $(X_1, T_1 \cup T_1', \emptyset)$ and $(X_2, \emptyset, T_2 \cup T_2')$, then the resulting amalgamation along $T$ has genus $3$ (see Figure~\ref{fig:amalgamationdiscrepancy} for a visual description of these two Heegaard splittings). There appears to be a discrepancy about which set of partitions for $X_1$ and $X_2$ to choose to obtain the amalgamation genus. 
\begin{figure}[h] 
   \centering
   \includegraphics[width=4.75in]{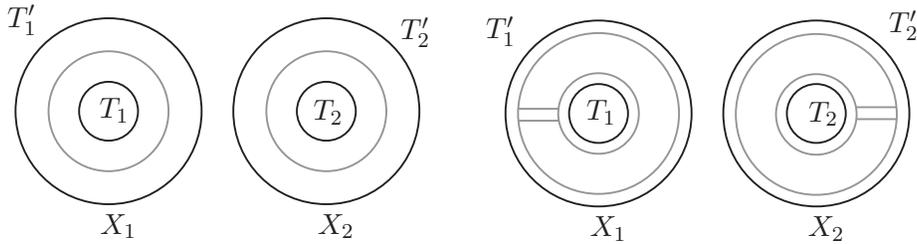} 
   \caption{Two amalgamations of $(M, T_1', T_2')$ along $T$.}
   \label{fig:amalgamationdiscrepancy}
\end{figure}

There are two ways of resolving this apparent ambiguity. If we are careful in the above example in keeping track of the compression bodies $V_i$ and $W_i$ for $X_i$, $i=1,2$, then it becomes clear that the first Heegaard splitting is actually a Heegaard splitting of $(M, T_2', T_1')$. Thus, if we regard $(M, \partial_1 M, \partial_2 M)$ and $(M, \partial_2 M, \partial_1 M)$ as different partitioned 3-manifolds, there is a unique partition of each which gives rise to the amalgamation genus. If, however, we consider $(M, \partial_1 M, \partial_2 M)$ and $(M, \partial_2 M, \partial_1 M)$ as the same partitioned 3-manifolds, then we have the choice of two sets of partitions for the components $X_i$, $1 \leq i \leq n$, as seen above. While both of these statements are equivalent, we will adopt the latter convention for the purposes of this paper. That is, we consider the partitioned 3-manifolds $(M, \partial_1 M, \partial_2 M)$ and $(M, \partial_2 M, \partial_1 M)$ as the same, and choose the set of partitions of $X_i$, $1 \leq i \leq n$, giving the smallest amalgamation genus. This convention is more natural topologically as the manifolds $(M, \partial_1 M, \partial_2 M)$ and $(M, \partial_2 M, \partial_1 M)$ are topologically the same, and also taking this convention works more smoothly for our arguments. 
\end{remark}

\begin{definition}
\rm
Of the two collections of partitions $(X_i, \partial_1 X_i, \partial_2 X_i)$, $1 \leq i \leq n$, discussed in Remark~\ref{twopartitionsremark}, the one which gives the amalgamation genus of $(M, \partial_1 M, \partial_2 M)$ along $F$ will be called {\em ideal}. If there is no discrepancy, then just choose one of the collections of partitions as ideal. 
\end{definition}

Let $F$ be an incompressible surface in $M$. Define a connected graph $\mathcal{G}_{M,F}$ associated to $(M,F)$ as follows. A vertex of $\mathcal{G}_{M,F}$ corresponds to a component of $M$ cut along $F$, and given two components $X_1$ and $X_2$ of $M$ cut along $F$, there is an edge  between their corresponding vertices for every component of $F$ contained in $\partial X_1 \cap \partial X_2$. Let $\ell = b_1(\mathcal{G}_{M,F}) = rank(H_1(\mathcal{G}_{M,F}; \mathbb{Z}))$. 

\begin{prop} 
\label{amalgamationgenusgeneral}
Let $F$ be a mutually separating incompressible surface in $M$, and let $X_1, \ldots, X_n$ be the components of $M$ cut along $F$. Then $$a(M, \partial_1 M, \partial_2 M, F) = \sum_{i=1}^n g_h(X_i, \partial_1 X_i, \partial_2 X_i) - g(F) + \ell $$
\noindent where $g(F)$ is the sum of the genera of each component of $F$, and the choice of the set of partitions $(X_i, \partial_1 X_i, \partial_2 X_i)$, $1 \leq i \leq n$, is ideal. 
\end{prop}

\begin{proof}
We induct on the number of components of $F$ amalgamated along. For some $1 \leq i, j \leq n$, $i \neq j$, let $M_0$ be the 3-manifold obtained by attaching $X_i$ to $X_j$ along a component $F_0$ of $F$ contained in $\partial X_i \cap \partial X_j$. Let $V_i \cup_{S_i} W_i$ and $V_j \cup_{S_j} W_j$ be minimal genus Heegaard splittings of $(X_i, \partial_1 X_i, \partial_2 X_i)$ and $(X_j, \partial_1 X_j, \partial_2 X_j)$, respectively. As described above, an amalgamation along $F_0$ is formed by compression bodies $V$ (obtained by attaching $1$-handles to $V_j$ and components of $(\partial_- V_i- F_0) \times I$) and $W$ (obtained by attaching $1$-handles to $W_i$ and components of $(\partial_- W_j- F_0) \times I$). A simple Euler characteristic argument gives that $a(M_0, \partial_1 M_0, \partial_2 M_0, F_0) = g_h(X_i, \partial_1 X_i, \partial_2 X_i) + g_h(X_j, \partial_1 X_j, \partial_2 X_j) - g(F)$.

Assume now that $m-1$ components of $F$ have been amalgamated along, and denote that subset of components as $F_{m-1}$. Let $M_{m-1}$ be the component of $M$ cut along $F - F_{m-1}$ containing $F_{m-1}$ (we may assume $M_{m-1}$ is connected). Then $M_{m-1}$ has a Heegaard splitting which is an amalgamation along $F_{m-1}$. Let $F_m$ be some component of $F - F_{m-1}$ which lies in $\partial M_{m-1}$, and let $X_i$ and $X_j$ be the components of $M$ cut along $F$ such that $F_m \subset \partial X_i \cap \partial X_j$. Without loss of generality, assume $X_i \subset M_{m-1}$ and set $M_m = M_{m-1} \cup X_j$ if $X_j$ is not contained in $M_{m-1}$, and $M_m = M_{m-1}$ otherwise. Then as above, we may form a Heegaard splitting of $(M_m, \partial_1 M_m, \partial_2 M_m)$ that is an amalgamation along $F_m \cup F_{m-1}$. 

The graph $\mathcal{G}_{M_m, F_m \cup F_{m-1}}$ is obtained from the graph $\mathcal{G}_{M_{m-1}, F_{m-1}}$ by adding an edge. If $b_1(\mathcal{G}_{M_m, F_m \cup F_{m-1}}) = b_1(\mathcal{G}_{M_{m-1}, F_{m-1}})$ after the addition of the edge, then $F_m$ cuts off $X_j$ from the rest of $M_m$. Hence the resulting Heegaard splitting of $M$ has genus $a(M_m, \partial_1 M_m, \partial_2 M_m, F_m \cup F_{m-1}) = a(M_{m-1}, \partial_1 M_{m-1}, \partial_2 M_{m-1}, F_{m-1}) + g_h(X_j, \partial_1 X_j, \partial_2 X_j) - g(F_m)$. If the addition of the edge results in $b_1(\mathcal{G}_{M_m, F_m \cup F_{m-1}}) = b_1(\mathcal{G}_{M_{m-1}, F_{m-1}}) + 1$, then both $X_i$ and $X_j$ are contained in $M_m$, and by an Euler characteristic argument, amalgamating along $F_m$ has the effect of decreasing the genus of the resulting Heegaard splitting by $g(F_m)$, but increasing it by $1$ due to the edge increasing the Betti number of the graph (as in what happens when computing an HNN extension of a fundamental group formed by taking a manifold and gluing together two of its boundary components). In other words, $a(M_m, \partial_1 M_m, \partial_2 M_m, F_m \cup F_{m-1}) = a(M_m, \partial_1 M_{m-1}, \partial_2 M_{m-1}, F_{m-1}) - g(F_m) + 1$. 

\begin{figure}[h] 
   \centering
   \includegraphics[width=5in]{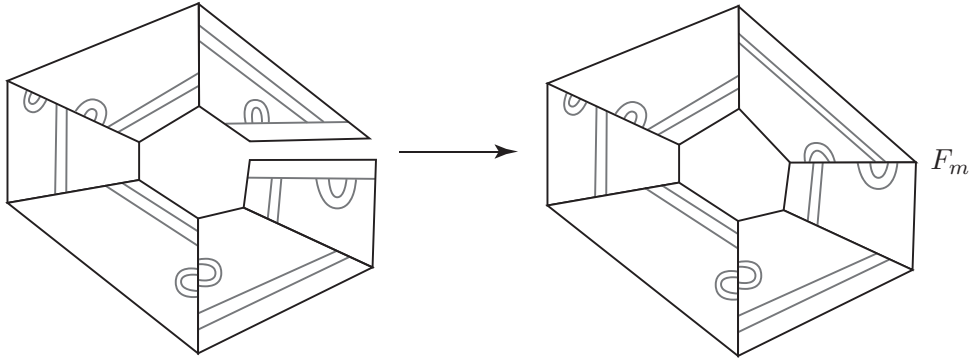} 
   \caption{Amalgamating along $F_m$ forms a cycle in $\mathcal{G}_{M_m, F_m}$ and changes $a(M_{m-1}, \partial_1 M_{m-1}, \partial_2 M_{m-1}, F_{m-1})$ by $-g(F_m) + 1$. }
   \label{fig:amalgamationgenus}
\end{figure}

\end{proof}

\subsection{Graph manifolds and amalgamations}

Heretofore our discussion has been for general 3-manifolds. We now shift attention to the case that $M$ is a totally orientable graph manifold with $\Theta$ a system of canonical tori in the JSJ decomposition. If $\mathcal{G}_{M, \Theta}$ has any cycles containing an odd number of vertices, then $\Theta$ fails to be mutually separating in $M$. 
Suppose that some cycle in $\mathcal{G}_{(M, \Theta)}$ contains an odd number of vertices. Choose some edge $e$ in this cycle, and take the corresponding torus $T_e$ and add a parallel copy $T_e'$ in $M$. Continue this process for additional cycles in $\mathcal{G}_{(M, \Theta)}$ containing an odd number of vertices. This yields a new disjoint union of tori $\Theta'$ which is mutually separating and contains $\Theta$ as a subset. Suppose that $\Theta'$ is obtained by adding $q$ such tori to $\Theta$, chosen to be minimal among all possible choices. 

\begin{definition}
\label{amalgamatablemodificationdefinition}
\rm The set $\Theta'$ constructed above will be called an {\em amalgamatable modification} of $\Theta$. 
\end{definition}

\begin{definition}
\label{amalgamationgenusdefinition2}
\rm For a partitioned graph manifold $(M, \partial_1 M, \partial_2 M)$, the {\em amalgamation genus}, denoted $a(M, \partial_1M, \partial_2 M)$, is defined to be $$a(M, \partial_1M, \partial_2 M) = a(M, \partial_1M, \partial_2 M, \Theta')$$ where $\Theta'$ is an amalgamatable modification of $\Theta$. 
\end{definition}

Note that if $M$ cut along $\Theta$ consists of a disjoint union of Seifert fibered components $X_1, \ldots, X_n$, then $M$ cut along $\Theta'$ consists of $X_1, \ldots, X_n$ (still referred to as the Seifert fibered components) as well as $q$ components which are homeomorphic to $T^2 \times I$. The graph $\mathcal{G}_{M,\Theta'}$ will be denoted simply as $\mathcal{G}_M$.

\begin{lemma}
\label{graph}
Let $\mathcal{G}_M$ be a connected graph with $n$ vertices and $m$ edges. Then $\ell = 1 - n + m$.
\end{lemma}

\begin{proof}
The Betti numbers of $\mathcal{G}_M$ are $b_0(\mathcal{G}_M) = 1$, $b_1(\mathcal{G}_M) = \ell$, and $b_j(\mathcal{G}_M) = 0$ for $j \geq 2$. The fact that $\chi(\mathcal{G}_M) = b_0(\mathcal{G}_M) - b_1(\mathcal{G}_M)$ implies that $n - m = 1 - \ell$, and hence $\ell = 1 - n +m$.
\end{proof}

Suppose $(M, \partial_1 M, \partial_2 M)$ is a partitioned graph manifold and suppose $M$ cut along $\Theta'$ has Seifert fibered components $X_1, \ldots, X_n$. Assume  the set of partitions $(X_i, \partial_1 X_i, \partial_2 X_i)$, $1 \leq i \leq n$, corresponds to the ideal set of partitions of components of $M$ cut along $\Theta'$. Define $a(X_i, \partial_1 X_i, \partial_2 X_i) = g_h (X_i, \partial_1 X_i, \partial_2 X_i) - 1$. 

\begin{lemma}
\label{graphamalgamationgenus}
If $(M, \partial_1 M, \partial_2 M)$ is a partitioned graph manifold with $\Theta'$ an amalagamatable modification of $\Theta$, then 

$$a(M, \partial_1 M, \partial_2 M) = \sum_{i=1}^n a(X_i, \partial_1 X_i, \partial_2 X_i) + q + 1.$$

\end{lemma}

\begin{proof}
Suppose that $T_1, \ldots, T_q$ are the tori in $\Theta$ that give parallel copies to form $\Theta'$. Let $p = |\Theta| - q$, so that $|\Theta'| = p + 2q$. Let $\ell = b_1(\mathcal{G}_M)$, and note that $g_h (T^2 \times I, T^2 \times \{ 0 \} \cup T^2 \times \{ 1 \}, \emptyset) = 2$. By Proposition~\ref{amalgamationgenusgeneral},  

$$a(M, \partial_1 M, \partial_2 M) = \sum_{i=1}^n g_h(X_i, \partial_1 X_i, \partial_2 X_i) + 2q - (p + 2q) + \ell$$

$$= \sum_{i=1}^n (a(X_i, \partial_1 X_i, \partial_2 X_i) + 1) - p + \ell$$

$$ = \sum_{i=1}^n a(X_i, \partial_1 X_i, \partial_2 X_i) + n - p + \ell.$$

By Lemma~\ref{graph}, $\ell = 1 - (n+ q) + (p + 2q) = 1 - n+ p + q$. Hence,

$$a(M, \partial_1 M, \partial_2 M) = \sum_{i=1}^n a(X_i, \partial_1 X_i, \partial_2 X_i) + q + 1.$$

\end{proof}

\subsection{Amalgamations and stabilization}
\label{amalgamationsandstabilization}

The following lemma establishes a relationship between stabilization and amalgamation for any compact, orientable 3-manifold.

\begin{lemma}
\label{stabilizationamalgamationlemma}
Let $F$ be a mutually separating incompressible surface in $M$, and let \hs \ and \pq \ be Heegaard splittings of $(M, \partial_1 M, \partial_2 M)$ which are both amalgamations along $F$. Suppose that $X_1, \ldots, X_n$ are the components of $M$ cut along $F$, and that $V_i \cup_{S_i} W_i$ and $P_i \cup_{\Sigma_i} Q_i$ are Heegaard splittings of $(X_i, \partial_1 X_i, \partial_2 X_i)$ obtained via untelescoping \hs \ and \pq , respectively. Let $\sigma_i$  be the minimal number of stabilizations of the higher genus splitting needed for $V_i \cup_{S_i} W_i$ and $P_i \cup_{\Sigma_i} Q_i$ to be isotopic in $X_i$. Then \hs \ and \pq \ are isotopic after at most $\max \{ \sigma_i \}$ stabilizations. 
\end{lemma}

\begin{proof}
Fix some $i$. Let $g_i$ be the genus of $V_i \cup_{S_i} W_i$ and $h_i$ the genus of $P_i \cup_{\Sigma_i} Q_i$. Without loss of generality, assume that $g_i \geq h_i$. Let $V_i' \cup_{S_i'} W_i'$ be the result of stabilizing $V_i \cup_{S_i} W_i$ $\sigma_i$ times, and $P_i' \cup_{\Sigma_i'} Q_i'$ be the result of stabilizing $P_i \cup_{\Sigma_i} Q_i$ $g_i - h_i + \sigma_i$ times. Let $\mathcal{B}$ be a 3-ball in $X_i$ such that $S_i' \cap \mathcal{B} = \Sigma_i' \cap \mathcal{B}$ is an equatorial disk in $\mathcal{B}$, and $V_i' \cap \mathcal{B} = P_i' \cap \mathcal{B}$. Then, the isotopy of $V_i' \cup_{S_i'} W_i'$ and $P_i' \cup_{\Sigma_i'} Q_i'$ taking $V_i'$ to $P_i'$, say, can be chosen to be fixed in $\mathcal{B}$. In fact, the same can be said for a finite number of 3-balls $\mathcal{B}_1, \ldots, \mathcal{B}_l$. 

Now, as \hs \ and \pq \ are amalgamations along $F$, we can assume that $S$ and $P$ intersect $F$ in $F$ minus disjoint open disks. Let $\mathcal{B}_1, \ldots, \mathcal{B}_l$ be pairwise disjoint 3-balls such that $\mathcal{B}_j \cap (F \cap \partial X_i)$ is an equatorial disk for $1 \leq j \leq l$, and each $\mathcal{B}_j$ contains an open disk in $(F \cap \partial X_i) - (S \cap P)$ where the 1-handle in $V$, $W$, $P$ or $Q$ meeting that open disk is not in $X_i$. Note that in $X_i - \cup_{j=1}^l \mathcal{B}_j$, \hs \ and \pq \ are equal to $V_i \cup_{S_i} W_i$ and $P_i \cup_{\Sigma_i} Q_i$ respectively. Thus, after $\sigma_i$ stabilizations, $V' \cap X_i$ is isotopic in $X_i$ to $P' \cap X_i$ by keeping $V' \cap \cup_{j=1}^l \mathcal{B}_j$ and $P' \cap \cup_{j=1}^l \mathcal{B}_j$ fixed. As \ppq \ has $\sigma_i$ stabilizations in $X_i$, those stabilizations can be isotoped into other components of $M$ cut along $F$ and the process can be repeated. Thus, \hs \ and \pq \ are isotopic after at most $\max \{ \sigma_i \}$ stabilizations. 

\end{proof}

\begin{corollary}
\label{graphmanifoldsamalgamationlemma}
Let \hs \ and \pq \ be Heegaard splittings of a partitioned graph manifold $(M, \partial_1 M, \partial_2 M)$ such that each splitting is an amalgamation along $\Theta'$, where $\Theta'$ is an amalgamatable modification of the system of canonical tori in the JSJ decomposition of $M$. Then \hs \ and \pq \ are isotopic after at most one stabilization of the higher genus splitting. 
\end{corollary}

\begin{proof}
Heegaard splittings of Seifert fibered spaces and manifolds homeomorphic to $T^2 \times I$ (with respect to appropriate partitions) are always isotopic after one stabilization (\cite{ScharlemannThompson}, \cite{Schultens1}). Hence, by Lemma~\ref{stabilizationamalgamationlemma}, \hs \ and \pq \ are isotopic after at most one stabilization of the higher genus splitting. 
\end{proof}

Corollary~\ref{graphmanifoldsamalgamationlemma} implies the case of Theorem~\ref{maintheorem} in which both \hs \ and \pq \ are amalgamations along $\Theta'$.

\subsection{Recognizing Amalgamations}

The following lemma gives a characterization of Heegaard splittings that are amalgamations along a mutually separating disjoint union of incompressible tori.

\begin{lemma}
\label{theamalgamationlemma}
Suppose $\Theta'$ is a mutually separating disjoint union of incompressible tori in a 3-manifold $M$ so that $M$ cut along $\Theta'$ is $M_1 \cup M_2$ (as in Definition~\ref{mutuallyseparatingdefinition}). Let \hs \ be a Heegaard splitting of $(M, \partial_1 M, \partial_2 M)$ such that $S$ is isotoped to meet each component of $\Theta'$ in curves which are essential in both surfaces (this can always be done {\em e.g.}~ for a stabilization of a strongly irreducible Heegaard splitting). Then $V \cup_S W$ is an amalgamation along $\Theta'$ if and only if (after isotopy) each annulus component of $V \cap \Theta'$ has a spanning disk in $V$ contained in $M_{\varepsilon}$, and each component of $W \cap \Theta'$ has a spanning disk in $W$ contained in $M_{\varepsilon'}$, where $\{ \varepsilon, \varepsilon' \} = \{1,2 \}$.
\end{lemma} 

The case of Lemma~\ref{theamalgamationlemma} where $\Theta'$ consists of a single torus is Lemma 3.1 in \cite{RDT1}. Applying that lemma to each component of $\Theta'$ gives the above result.

\section{Heegaard splittings of graph manifolds and the active component}
\label{sec:section4}

For the remainder of this paper, $M$ will denote a totally orientable graph manifold. For this section, consider $M$ as discussed in Remark~\ref{graphmanifoldsremark} with $N(\Theta)$ as edge manifolds and $\overline{M - N(\Theta)}$ as Seifert fibered components. We refer to an edge manifold as $N(T)$ where $T$ is a component of $\Theta$.

\begin{remark}
\rm
It follows from Theorem 4.1 in \cite{BSS} that generically, Heegaard splittings of graph manifolds are amalgamations along incompressible surfaces. In a graph manifold, a component of an incompressible surface can be isotoped to be a canonical torus, or to be vertical or horizontal in each Seifert component of the manifold (see Lemma 4.2 in \cite{Schultens2}). The latter situation is somewhat restricted by the gluings of the Seifert fibered components along the canonical tori, unless the incompressible surface is a union of tori (possibly canonical) that are vertical in individual Seifert fibered components. If this is the case, cutting along these tori yields graph manifolds containing Heegaard splitting obtained from untelescoping, which can be analyzed in an inductive fashion. Thus, it follows that generically, Heegaard splittings of graph manifolds are amalgamations along an amalgamatable modification of the system of canonical tori in the JSJ decomposition. By Corollary~\ref{graphmanifoldsamalgamationlemma} any two such splittings are isotopic after at most one stabilization. Thus, the Heegaard splittings which are not amalgamations along an amalgamatable modification of the system of canonical tori in the JSJ decomposition are of interest.
\end{remark}

\begin{definition}
\label{vhpvph}
\rm A (not necessarily connected) surface $F$ in a Seifert fibered space $X$ is {\em vertical} if it is incompressible and consists of fibers of $X$, and is {\em horizontal} if it is incompressible and transverse to each fiber of $X$. A surface $F$ in $X$ is {\em pseudovertical} if it is obtained from a vertical surface via ambient 1-surgery along arcs which meet the vertical surface at their endpoints and which project to distinct simple arcs in the base $B$ of $X$.  If $f$ is a fiber of $X$, $F$ is {\em pseudohorizontal} if it is horizontal in $\overline{X - N(f)}$ and $F \cap N(f)$ is an incompressible annulus in $N(f)$ with core $f$.
\end{definition}

The main theorem in \cite{Schultens2} is that any strongly irreducible Heegaard splitting \hs \ of $M$ has a standard form.

\begin{theorem}[\cite{Schultens2}, Theorem 1.1]
\label{Schultenstheorem}
Let \hs \ be a strongly irreducible Heegaard splitting of a totally orientable partitioned graph manifold $(M, \partial_1 M, \partial_2 M)$. Then \hs \ can be isotoped to intersect each Seifert fibered component $X$ of $M$ such that $S \cap X$ is horizontal, pseudohorizontal, vertical, or pseudovertical, and $S$ intersects each edge manifold $N(T)$ such that:

\begin{itemize}
\item[1.] $S \cap N(T)$ is obtained from a disjoint union of incompressible annuli in $N(T)$ by possibly performing ambient 1-surgery along an arc isotopic into $\partial N(T)$, or

\item[2.] There are a pair of simple closed curves $c_0$ and $c_1$ on $T$ that intersect in a single point $p$, and $S \cap N(T)$ is isotopic to the frontier of a regular neighborhood of $c_0 \times \{0\} \cup \{p\} \times I \cup c_1 \times \{1\}$.
\end{itemize}
\end{theorem}

In the case that $\Theta$ is empty, the above theorem reduces to the characterization of Heegaard splittings of Seifert fibered spaces given in \cite{MoriahSchultens}. As the stabilization problem has been solved for Heegaard splittings of Seifert fibered spaces in \cite{Schultens1}, we shall assume henceforth that $\Theta$ is nonempty.

Let \hs \ be a strongly irreducible Heegaard splitting of $(M, \partial_1 M, \partial_2 M)$. In obtaining the above characterization of \hs \ it is useful to isotope $S$ so that $S$ is essential in all but one component of $\overline{M - N(\Theta)} \sqcup N(\Theta)$. 

\begin{lemma}[\cite{Schultens2}, Lemma 6.1]
\label{activecomponentlemma}
Let \hs \ be a strongly irreducible Heegaard splitting of a partitioned graph manifold $(M, \partial_1 M, \partial_2 M)$. If $D \subset V$ and $E \subset W$ are essential disks, then, after isotopy of $S$, the outermost disk components of $D$ cut along $\Theta$ and $E$ cut along $\Theta$ lie in the same Seifert fibered component or edge manifold $N$ of $M$. Moreover, for each Seifert fibered component or edge manifold $N' \neq N$, $S \cap N'$ is incompressible.
\end{lemma}

The proof relies on the fact that $\Theta$ consists of tori. The idea is that by taking an outermost disk component $D'$ of $D$ cut along $\Theta$ and the annulus $A$ of $V \cap \Theta$ which it meets, two copies of $D'$ attached to a rectangle in $A$ can be isotoped to a compressing disk for $V$ inside the Seifert fibered piece or edge manifold $N$ containing $D'$. The fact that it is a compressing disk for $V$ in $N$ follows from the assumption that $|S \cap \Theta|$ is minimal (see Remark~\ref{startingpositionremark}) so that $S \cap N$ contains no boundary parallel annuli. The strong irreducibility of \hs \ then implies that all such disks must exist in $N$, implying the lemma.

\begin{definition}
\label{activecomponentdefinition}
\rm The component $N$ in Lemma~\ref{activecomponentlemma} is called the {\em active component}. 
\end{definition}

Further examination of the active component reveals that its structure can be well understood.

\begin{lemma}
\label{edgemanifoldrefinement1}
Let \hs \ be a strongly irreducible Heegaard splitting of a partitioned graph manifold $(M, \partial_1 M, \partial_2 M)$. After isotopy, either the active component $N$ is an edge manifold, or it is a Seifert fibered component and $S \cap N$ is psuedohorizontal.
\end{lemma}

\begin{proof}
This lemma follows straightforwardly from Proposition 7.5 and Lemma 7.10 in \cite{Schultens2}.
\end{proof}

Upon further inspection, we have an explicit understanding of what the Heegaard splitting surface looks like in the active component if it is an edge manifold.

\begin{lemma}
\label{activecomponentrefinement2}
Let \hs \ be a strongly irreducible Heegaard splitting of $(M, \partial_1 M, \partial_2 M)$, and suppose that the active component $N$ is an edge manifold. Then after isotopy, either:

\begin{itemize}
\item[1.] $S \cap N$ is a disjoint union of essential annuli and a component obtained by taking two incompressible annuli, at least one of which is boundary parallel in $N$, and performing ambient 1-surgery along an arc parallel into the boundary connecting the two annuli, or

\item[2.] There are a pair of simple closed curves $c_0$ and $c_1$ on $T$ that intersect in a single point $p$, and $S \cap N$ is isotopic to the frontier of a regular neighborhood of $c_0 \times \{0\} \cup \{p\} \times I \cup c_1 \times \{1\}$. 
\end{itemize}
\end{lemma}

\begin{proof}
If \hs \ is a strongly irreducible Heegaard splitting of $(M, \partial_1 M, \partial_2 M)$, then as in Remark~\ref{startingpositionremark}, $S$ can be isotoped so that no component of $V \cap \Theta$ or $W \cap \Theta$ (and hence no component of $ V \cap \partial N(\Theta)$ or $W \cap  \partial N(\Theta)$) is parallel into $\partial_+ V$ or $\partial_+ W$. Furthermore, the isotopies performed in Proposition 7.5 and Lemma 7.10 in \cite{Schultens2} to give Lemma~\ref{edgemanifoldrefinement1} do not introduce any boundary parallel annuli into the active component when it is an edge manifold. Thus, we can assume that $S \cap N$ does not contain any components which are parallel into $\partial N$. 
The only remaining possibility given by Theorem~\ref{Schultenstheorem} is that $S \cap N$ consists of a disjoint union of essential annuli in $N$ and a 4-times punctured sphere component, obtained from ambient 1-surgery along an arc $a$ connecting two incompressible annuli in $N$. 

Suppose $D \subset V$ and $E \subset W$ are essential disks. Then $D \cap \partial N$ and $E \cap \partial N$ are nonempty, and by Lemma~\ref{activecomponentlemma} the outermost disk components of $D$ and $E$ cut along $\partial$ lie in $N$ as boundary compressing disks of $S \cap N$. If $a$ connects two essential annuli in $N$, then there are boundary compressing disks for $S \cap N$ in one of the compression bodies but not the other (see Figure~\ref{fig:activecomponentrefinement}), contradicting the existence of boundary compressing disks given by outermost disk components of both $D$ and $E$. Thus we conclude that at least one of the annuli meeting $a$ is boundary parallel. 

\begin{figure}[h] 
   \centering
   \includegraphics[width=1.5in]{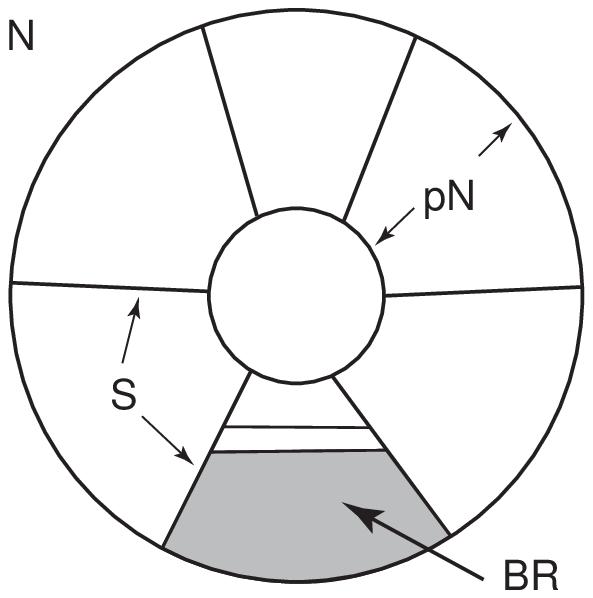} 
   \caption{$S$ cannot have a component obtained from ambient 1-surgery along two essential annuli in $N$.}
   \label{fig:activecomponentrefinement}
\end{figure}

Note that if the annuli meeting $a$ are both boundary parallel into different components of $\partial N$ such that the projection of their cores on a component of $\partial N$ has the same slope, then the 4-punctured sphere component of $S \cap N$ is isotopic to the 4-punctured sphere component obtained from performing ambient 1-surgery along an arc connecting two essential annuli in $N$ (see {\em e.g.}~Figure~\ref{fig:dualtubeswap}), which is the case above.  

\end{proof}

We summarize the above results in the following definition.

\begin{definition}
\label{activeposition}
\rm Suppose that \hs \ is a strongly irreducible Heegaard splitting of a paritioned graph manifold $(M, \partial_1 M, \partial_2 M)$ isotoped as in Lemma~\ref{edgemanifoldrefinement1} and Lemma~\ref{activecomponentrefinement2}. Then \hs \ is said to be in {\em active position} with respect to $\Theta$. In particular, $S$ is essential in each component of $(\overline{M - N(\Theta)} \sqcup N(\Theta)) - N$, and $N$ is either an edge manifold, or $N$ is a Seifert fibered component and $S \cap N$ is pseudohorizontal. 
\end{definition}

\begin{figure}[h] 
   \centering
   \includegraphics[width=2.5in]{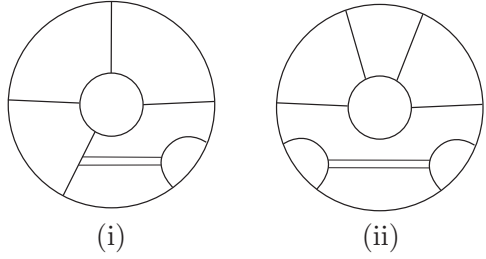} 
   \caption{Possibilities (i) and (ii) (schematic) for conclusion 1 in Lemma~\ref{activecomponentrefinement2} ($S \cap N$ is aligned).}
   \label{fig:option1}
\end{figure}

\begin{figure}[h] 
   \centering
   \includegraphics[width=1.2in]{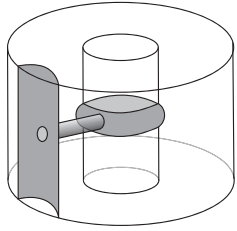} 
   \caption{Possibility for conclusion 2 in Lemma~\ref{activecomponentrefinement2} ($S \cap N$ is a  toggle).}
   \label{fig:option2}
\end{figure}

If $N$ is an edge manifold, then either $S \cap N$ is {\em aligned} (Figure~\ref{fig:option1}), or $S \cap N$ is a {\em toggle} (Figure~\ref{fig:option2}).

\bigskip

Suppose that \hs \ is in active position with respect to $\Theta$. Then the active component $N$ is either equal to $N(T)$ for some component $T$ of $\Theta$, or $N$ is a Seifert fibered component of $M$ such that $S \cap N$ is psuedohorizontal. In the former case, to simplify our arguments, assume that $\Theta$ is modified by replacing $T$ with $\partial N(T)$. This allows us to assume that \hs \ being in active position implies for every component $X$ of $M$ cut along $\Theta$ except for $N$, $S \cap X$ is either vertical or horizontal in the context of Definition~\ref{vhpvph}. 
\begin{definition}
\label{verticalhorizontalcomponents}
\rm Suppose that, if necessary, $\Theta$ is modified as above. A component $X$ of $M$ cut along $\Theta$ is a {\em horizontal component} if $S \cap X$ is horizontal, and a {\em vertical component} if $S \cap X$ is vertical. 
\end{definition}

Note that this definition only makes sense when $S$ is in a fixed active position. Also note that components of $\partial M$ can only lie in vertical components. 

\section{The case of no horizontal components} 
\label{sec:section5}

The goal of the next sections will be to prove the following generalization of Theorem~\ref{maintheorem}.

\begin{theorem}
\label{maintheorem1}
Let $(M, \partial_1 M, \partial_2 M)$ be a partitioned totally orientable graph manifold with $\Theta$ a system of canonical tori in the JSJ decomposition, and let $\Theta'$ be an amalgamatable modification of $\Theta$. Suppose \hs \ and \pq \ are Heegaard splittings of $(M, \partial_1 M, \partial_2 M)$ such that \hs \ either is an amalgamation along $\Theta'$ or is strongly irreducible, and \pq \ is an amalgamation along $\Theta'$. Then \hs \ and \pq \ are isotopic after at most one stabilization of the larger genus splitting.
\end{theorem}

As mentioned above, we assume that $\Theta' \neq \emptyset$. The case that \hs \ and \pq \ are both amalgamations along $\Theta'$ is established in Corollary~\ref{graphmanifoldsamalgamationlemma}. Thus assume that \hs \ is a strongly irreducible Heegaard splitting of genus $g$. As $\Theta'$ is mutually separating, denote $M_1$ and $M_2$ as the manifolds obtained from cutting $M$ along $\Theta'$ (see Definition~\ref{mutuallyseparatingdefinition}). We will use the structure of \hs \ provided in Section~\ref{sec:section4} to show that if $g \geq a(M, \partial_1M, \partial_2 M)$, then after stabilizing \hs \ once to \hhs, $S'$ can be isotoped so that every component of $V' \cap \Theta'$ has a spanning disk in $M_{\varepsilon}$ and every component of $W' \cap \Theta'$ has a spanning disk in $M_{\varepsilon'}$, where $\{ \varepsilon, \varepsilon ' \} = \{ 1, 2\}$. Lemma~\ref{theamalgamationlemma} then implies that \hhs \ is an amalgamation along $\Theta'$. If $g < a(M, \partial_1 M, \partial_2 M)$, then the same result holds after stabilizing \hs \ $a(M, \partial_1 M, \partial_2 M) - g + 1$ times. 
\smallskip

We use this strategy first to prove Theorem~\ref{maintheorem1} in the special case that $M$ cut along $\Theta$ has no horizontal components. 

\begin{lemma}
\label{nohorizontalcomponentslemma}
If $M$ cut along $\Theta$ has no horizontal components and if $N$ is an edge manifold (in particular $S\cap N$ is not psuedohorizontal), then \hs \ and a Heegaard splitting which is an amalgamation along $\Theta'$ are isotopic after at most one stabilization of the higher genus splitting. 
\end{lemma}

\begin{proof} By assumption, $M$ cut along $\Theta$ contains only vertical components, and the active component $N$ is an edge manifold. Note that for two vertical components $Y_1$ and $Y_2$, it must be the case that $\partial Y_1 \cap \partial Y_2 = \emptyset$. Otherwise, the fiberings of $Y_1$ and $Y_2$ would line up and they would be considered a single vertical component of $M$. Thus, other than the genus 1 splitting of $T^2 \times I$ (the only example of a strongly irreducible amalgamation along a torus), there are only two possible configurations for $M$: $(1)$ either $M = Y_1 \cup N \cup Y_2$ where $N$ is an edge manifold and $S \cap N$ is a toggle (see Figure \ref{fig:option2}), $Y_i$ is a vertical component and $S \cap Y_i$ is a single vertical annulus for $i = 1,2$, or $(2)$ $M = Y_1 \cup N$ where $S \cap N$ is a toggle and $Y_1$ is a vertical component with $S \cap Y_1$ a disjoint union of two vertical annuli, each having a boundary component on each component of $\partial Y_1$. Note that this implies the genus of \hs \ is 2. In case $(1)$, for each $i=1,2$, $Y_i$ is a Seifert fibered space over a disk with two exceptional fibers, over an annulus with one or zero exceptional fibers, or over a twice punctured disk with no exceptional fibers. In case $(2)$, $Y_1$ is a Seifert fibered space over an annulus with zero, one or two exceptional fibers, over a twice punctured disk with one or zero exceptional fibers, or over a three-times punctured disk with no exceptional fibers. 

For $(1)$, one can check that $1 \leq a(M, \partial_1 M, \partial_2 M) \leq 3$, depending on the partition $(M, \partial_1 M, \partial_2 M)$. Consider $\Theta$ as one of the boundary components of $N$, without loss of generality the component of $N \cap Y_1$. The compression bodies $V$ and $W$ intersect $\Theta$ in two annuli $A_V$ and $A_W$, respectively. Note that $A_W$ has a spanning disk in $N$ given by the toggle. Stabilize \hs \ by adding a 1-handle $H$ to $W$ in $Y_1$ such that $H$ is adjacent to $A_V$. Then there is a spanning disk of $A_V$ which intersects the cocore of $H$ in a single point. By Lemma~\ref{theamalgamationlemma}, the stabilized splitting \hhs \ is an amalgamation along $\Theta$.

\begin{figure}[h] 
   \centering
   \includegraphics[width=1in]{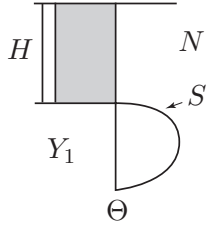} 
   \caption{An amalgamation along $\Theta$.}
   \label{fig:twoverticalcomponents}
\end{figure}

For $(2)$, $\Theta$ is non-separating, so we can take $\Theta' = \partial N$. Here, since $\partial Y_1$ contains the two components of $\Theta'$, it is straightforward to check that $3 \leq a(M, \partial_1 M, \partial_2 M) \leq 4$, with $a(M, \partial_1 M, \partial_2 M) = 3$ if and only if $Y_1$ is a Seifert fibered space over an annulus with zero or one exceptional fibers, or over a twice punctured disk with no exceptional fibers. In these cases, stabilizing $S$ in $Y_1$ as above provides spanning disks in $Y_1$ for two of the annuli in $V \cap \Theta'$ (see Figure~\ref{fig:onestabilizationtwodisks}). The toggle gives spanning disks in $N$ for the components of $W \cap \Theta'$. By Lemma~\ref{theamalgamationlemma}, \hhs \ is isotopic to an amalgamation along $\Theta'$.

\begin{figure}[h] 
   \centering
   \includegraphics[width=2.5in]{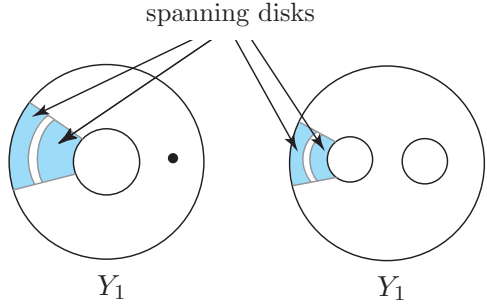} 
   \caption{One stabilization of $S$ in $Y_1$ gives two spanning disks for components of $V \cap \Theta'$.}
   \label{fig:onestabilizationtwodisks}
\end{figure}

In the case that $a(M, \partial_1 M, \partial_2 M) = 4$, then after 2 stabilizations, one for each component of $\Theta'$ as done for case $(1)$, $V' \cup_{S'} W'$ is an amalgamation along $\Theta'$ by Lemma~\ref{theamalgamationlemma}. 

In both cases $(1)$ and $(2)$, the Heegaard splitting \hs \ becomes an amalgamation along $\Theta'$ after stabilizing the splitting to have genus $a(M, \partial_1 M, \partial_2 M)$. By Corollary~\ref{graphmanifoldsamalgamationlemma}, any two amalgamations along $\Theta'$ are isotopic after at most 1 stabilization. 

\end{proof}

\section{Horizontal components and tubes}
\label{sec:horizontalcomponents}

Having established Theorem~\ref{maintheorem1} in the case that $M$ has no horizontal components and the active component is an edge manifold, assume henceforth that $M$ has a horizontal component or that the active component $N$ is Seifert fibered and $S \cap N$ is pseudohorizontal. 

Let $\Theta'$ be an amalgamatable modification of $\Theta$, constructed as in Defintion~\ref{amalgamatablemodificationdefinition}. Note that $M$ cut along $\Theta'$ consists of Seifert fibered components homeomorphic to the components of $M$ cut along $\Theta$, in addition to components that are homeomorphic to $T^2 \times I$. As before, we will refer to the Seifert fibered components as horizontal, vertical or active, and the remaining components (other than the active component) as {\em $T^2 \times I$ components}. Note that by definition, no component of $M$ cut along $\Theta'$ meets a component of $\Theta'$ on both sides. 

\subsection{Horizontal surfaces and fundamental 2-complexes}

Let $X$ be a horizontal component of $M$, and suppose $X$ has Seifert fibered data given by $$\left( g; m; \frac{\beta_1}{\alpha_1}, \ldots, \frac{\beta_k}{\alpha_k} \right).$$

\noindent Label the boundary components of the base $B$ of $X$ by $\delta_1, \ldots, \delta_m$. Let $\lambda_1, \ldots, \lambda_{2g}$ be properly embedded arcs in $B$ from $\delta_1$ to itself, $\omega_2, \ldots, \omega_m$ properly embedded arcs in $B$ from $\delta_1$ to $ \delta_2, \ldots, \delta_m$ respectively, chosen so that the arcs are pairwise disjoint and cut $B$ into a disk. If $k \geq 2$ let $\gamma_1, \ldots, \gamma_{k-1}$ be arcs in $B$ from $\delta_1$ to $\delta_1$ which cut off each of the exceptional points $e_1, \ldots , e_{k-1}$ from the rest of $B$, disjoint from the previous arcs.
\begin{figure}[h] 
   \centering
   \includegraphics[width=2 in]{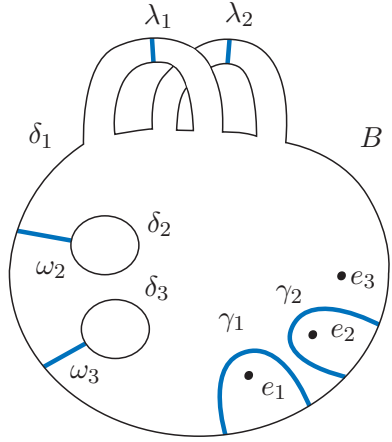} 
   \caption{Arcs that cut $B$ into $k-1$ disks each containing an exceptional point if $k \geq 2$, or one disk with one or no exceptional points if $k \leq 1$.}
   \label{fig:basespanningarcs}
\end{figure}
Let $T_1$ denote the component of $\partial X$ which projects to $\delta_1$ on $B$. 

\begin{definition}
\rm Each of the arcs $\lambda_1, \ldots, \lambda_{2g}, \omega_2, \ldots, \omega_m, \gamma_1, \ldots, \gamma_{k-1}$ lifts to an essential annulus $\Lambda_i = \lambda_i \times S^1$, $\Omega_i =  \omega_i \times S^1$, and $\Gamma_i = \gamma_i \times S^1$ in $X$. The union of these annuli 

\begin{displaymath}
G = \bigcup_{i=1}^{2g}\Lambda_i \cup \bigcup_{i=2}^m \Omega_i \cup \bigcup_{i=1}^{k-1} \Gamma_i
\end{displaymath}

\noindent is a called a {\em fundamental 2-complex} for $X$ built at $T_1$. 
\end{definition}

The main feature of a fundamental 2-complex $G$ for $X$ is that $X$ cut along $G$ is a disjoint union of solid tori with an exceptional fiber in each solid torus if $k \geq 1$, or a single solid torus containing no exceptional fibers if $k=0$. In particular, if $F$ is a connected horizontal essential surface in $X$, $F$ cut along the arcs in $G \cap F$ is a disjoint union of disks. 

\subsection{Tube sliding} The goal of this section, established in Lemma~\ref{tubemigrationlemma}, is to show that stabilizing \hs \ in a horizontal component produces tubes which can be exported into neighboring components of $M$ cut along $\Theta'$. The following discussion describes the types of isotopy that we use to show this result. Note that for convenience we will often fix a compression body in making our arguments ({\em e.g.}~stabilizing in $V$), although a symmetric result holds when choosing the other compression body instead.
\smallskip

Let $X$ be a horizontal component of $M$ cut along $\Theta'$, and let $Z_1, \ldots, Z_n$ be the components of $V \cap X$. Then for each $1 \leq i \leq n$, $Z_i$ is homeomorphic to $F \times I$, where $F$ is a punctured surface (in particular $Z_i$ is a handlebody). Fix some $i$. By Proposition~\ref{spanningdisk} every component of $Z_i \cap \partial X$ has a spanning disk in $V$. Then, taking outermost disk components if necessary, some component $A$ of $Z_i \cap \partial X$ has a spanning disk $D$ such that $D \cap Z_i$ is a spanning arc $a$ of $A$. 

Let $\alpha$ be a vertical arc in $Z_i$ parallel to $a$. Then the disk $D$ gives that $\alpha$ is boundary parallel. Thus, we can stabilize \hs \ in $Z_i$ by attaching a 1-handle $H$ to $W$ such that the core of $H$ is $\alpha$ and removing $H$ from $V$. Let \hhs \ denote this stabilization and note that $S = S'$ outside of $Z_i$, while $S' \cap Z_i$ is connected. 

\begin{figure}[h] 
   \centering
   \includegraphics[width=5in]{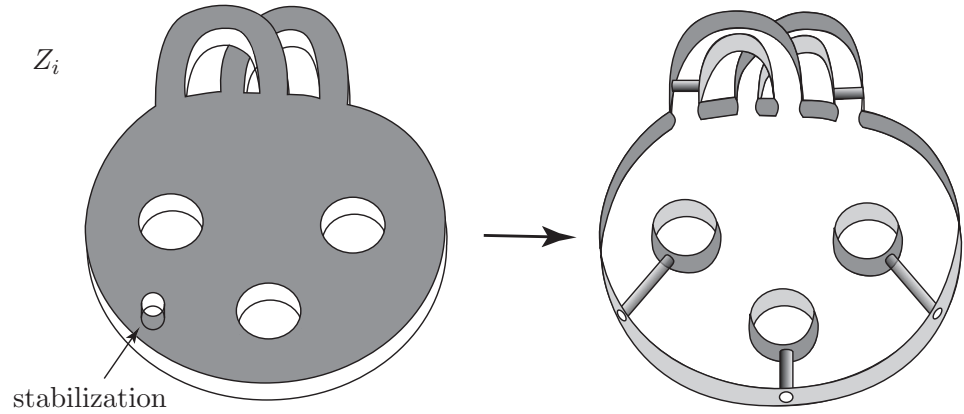} 
   \caption{Isotoping $S'$ in $Z_i$.}
   \label{fig:horizontalisotopy}
\end{figure}

Observe that $S' \cap Z_i$ is isotopic $(\mbox{rel} \ \partial)$ to the surface obtained from a disjoint union of properly embedded annuli in $Z_i$ parallel to the components of $Z_i \cap \partial X$ via ambient 1-surgery along arcs connecting those annuli (see Figure~\ref{fig:horizontalisotopy}). Let $\mathcal{T}$ denote the arcs used for the ambient 1-surgery, and take $\mathcal{T}$ to be properly embedded in $ F_0 = F \times \{1/2\}$, where $Z_i = F \times I$. Note that $F_0$ cut along $\mathcal{T}$ is a disk. 

Let $G$ be a fundamental 2-complex of $X$ built at $T_1$, where $T_1$ is some component of $\partial X$. As $G$ cuts $F_0$ into disks, the arcs in $\mathcal{T}$ can be chosen to be a subset of the arcs in $G \cap F_0$; namely all of the arcs in $\Lambda \cap F_0$, all of the arcs in $\Omega \cap F_0$, and a subset of the arcs in $\Gamma \cap F_0$ such that the union of the chosen arcs does not separate $F_0$. 

\begin{definition}
\rm Let $D$ be a compressing disk for $\partial_+ V'$ in $V'$ (of $\partial_+ W'$ in $W'$). A {\em tube} in $V'$ (in $W'$) is a regular neighborhood $D \times I$ of $D$ such that $\partial D \times I \subset \partial_+ V'$ ($\subset \partial_+ W'$). For a tube $\tau = D \times I$, the {\em core} of $\tau$ is an arc $p_0 \times I$ where $p_0$ is a point chosen in the interior of $D$, the {\em cocore} of $\tau$ is the simple closed curve $\partial D \times \{1/2\}$, and the {\em feet} of $\tau$ are the disks $D \times \{0\}$ and $D \times \{1\}$.
\end{definition}

\begin{definition}
\rm After performing the above isotopy of $S' \cap Z_i$ so that the arcs in $\mathcal{T}$ are spanning arcs of components of $G$ and form cores of tubes in $V'$ (as in Figure~\ref{fig:horizontalisotopy}), we say that $S' \cap Z_i$ is in {\em tube position} with respect to $G$.
\end{definition}

\begin{figure}[h]
   \centering
   \includegraphics[width=1.75in]{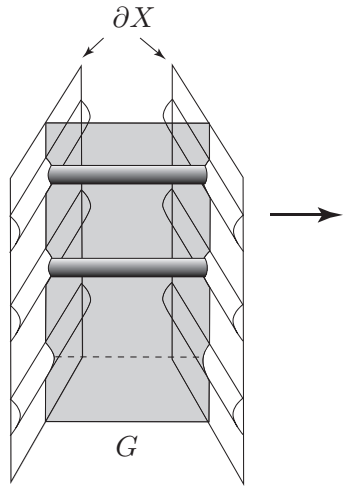} 
   \caption{The tubes have cores on $G$.}
   \label{fig:tubesoncomplex}
\end{figure}

Note that above we obtained $S' \cap Z_i$ in tube position by stabilizing $V$ in $Z_i$. The same positioning of $S' \cap Z_i$ could have been obtained if some tube in $W'$ were isotoped into $Z_i$ so that its core is parallel to a spanning arc of a component of $Z_i \cap \partial X$. 

\begin{remark}
\label{pseudohorizontalremark}
\rm Assume that $X = N$ is a Seifert fibered component of $M$ cut along $\Theta'$ which is the active component. By Lemma~\ref{edgemanifoldrefinement1}, $S \cap N$ is pseudohorizontal. Thus $S \cap N$ is connected and $S \cap (\overline{N - N(f)})$ is two copies of a horizontal surface in $\overline{N - N(f)}$, where $f$ is some fiber in $N$. Isotoping the annulus $S \cap N(f)$ into $\overline{N - N(f)}$ allows $S \cap N$ to be isotoped in $N$ to be in tube position (without adding a tube), just as was the case for a horizontal component above (after adding a tube). This allows us to treat $N$ as a horizontal component by considering the structure of $S \cap N$ in tube position.
\end{remark}

\begin{definition}
\rm
Two tubes $\tau_1$ and $\tau_2$ are {\em adjacent} on $G$ if the cores of $\tau_1$ and $\tau_2$ lie on the same annulus component of $G$, and there are no arcs of $G \cap F_0$ between them.
\end{definition}

The arcs of $G \cap F_0$ used to form $\mathcal{T}$ can be chosen so that all of the arcs of a particular $\Gamma_i \cap F_0$ are adjacent on $\Gamma_i$. As $\mathcal{T}$ contains all the arcs in $\Lambda \cap F_0$ and $\Omega \cap F_0$, we can assume that all arcs in $\mathcal{T}$ on some component of $G$ are adjacent. 
\smallskip

For the following definition, first assume that $V \cap X$ has only one component $Z_1$.

\begin{figure}[h] 
   \centering
   \includegraphics[width=4.5in]{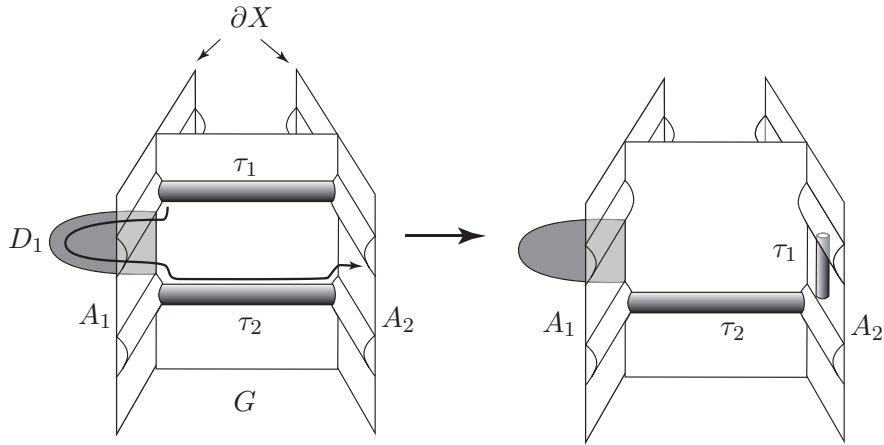} 
   \caption{A tube slide.}
   \label{fig:tubeslide}
\end{figure}

\begin{definition}
\label{tubeslidedefinition}
\rm Suppose that $\tau_1$ and $\tau_2$ are adjacent tubes on $G$. Let $A_1$ and $A_2$ be annulus components of $W' \cap \partial X$ between the feet of $\tau_1$ and $\tau_2$ (possibly $A_1 = A_2$). If there exists a spanning disk $D_1$ for $A_1$ such that $D_1 \cap X$ is a spanning arc of $A_1$, then we may isotope $\tau_1$ along a neighbrohood of $D_1$ and a rectangle in $G - \tau_1 \cup \tau_2$ so that the core of $\tau_1$ becomes a spanning arc of $A_2$. Such a move is called a {\em tube slide} into $A_2$ (see Figure~\ref{fig:tubeslide}).
\end{definition}

Let $T$ be the component of $\partial X$ containing $A_2$. After slightly pushing $\tau_2$ through $A_2$, the essential curves of intersection $S' \cap X$ are the same as before the tube slide. Hence after tube sliding, $V' \cap \Theta = V \cap \Theta$ and $W' \cap \Theta = W \cap \Theta$ remain a disjoint union of properly embedded incompressible annuli, and $\tau_2$ is contained in the component of $M$ cut along $\Theta'$ across $T$ from $X$. 
\smallskip

\begin{figure}[h] 
   \centering
   \includegraphics[width=4.5in]{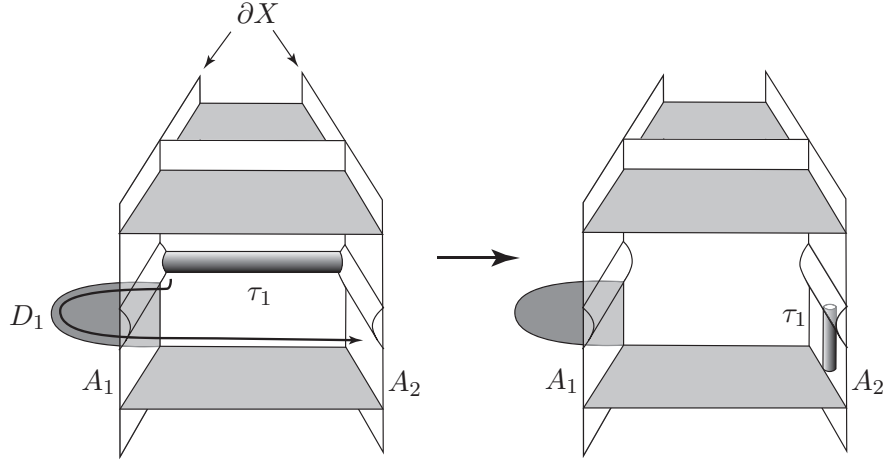} 
   \caption{A tube slide in the case $V \cap X$ has more than 1 component.}
   \label{fig:tubeslide2}
\end{figure}

Suppose now that $V \cap X$ has more than one component. Let $Z_i$ be such a component, and stabilize \hs \ in $Z_i$ as before. Set $\hat{S'} = (S' \cap X) -  (S' \cap Z_i)$. Note that $\hat{S'}$ is a disjoint union of parallel horizontal surfaces in $X$. After isotoping $Z_i \cap S'$ into tube position, as above we may perform a tube slide for any of the resulting tubes. This time, however, the tube $\tau_1$ is isotoped through a rectangle in $G - (\tau_1 \cup \hat{S'})$ and a foot of $\tau_1$ is isotoped along a component of $\hat{S'}$ (see Figure~\ref{fig:tubeslide2}).

\begin{definition}
\rm Let $\tau$ be a tube in $V'$, say, and suppose that $\tau$ is contained in a component $X$ of $M$ cut along $\Theta'$. Let $A$ be a component of $W' \cap \partial X$. Then $\tau$ is {\em adjacent to $A$} if there as a spanning disk of $A$ contained in $W' \cap X$ meeting the cocore of $\tau$ in a single point.  
\end{definition}

Note that $\tau$ being adjacent to $A$ is equivalent to the condition that $\tau$ can be isotoped so that its core is a spanning arc of $A$, and this isotopy is performed so that the core of $\tau$ is always in $X$. 

\begin{figure}[h] 
   \centering
   \includegraphics[width=.6in]{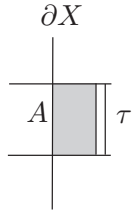} 
   \caption{A tube $\tau$ adjacent to the annulus $A$.}
   \label{fig:tubespanningdisk2}
\end{figure}

The following discussion describes a second type of isotopy we will use for isotoping tubes. 

\begin{definition}
\rm A tube $\tau$ in $V'$ is {\em trivial} if it can be isotoped so that there is an essential disk in $W'$ meeting the cocore of $\tau$ in $V'$ in a single point. 
\end{definition}

\begin{lemma}
\label{trivialtubelemma}
Suppose that $\tau$ is a trivial tube in either $V'$ or $W'$. Then $\tau$ can be isotoped so that it is adjacent to any component of $V' \cap \Theta$ or any component of $W' \cap \Theta$.
\end{lemma}

\begin{proof}
The proof is similar to showing that there is only one way to stabilize a given Heegaard splitting up to isotopy. For a fixed compressing disk in $W'$, say, a trivial tube in $V'$ can be isotoped such that its core becomes a properly embedded arc in this disk. Let $\tau$ be a trivial tube in $V' \cap X$, and let $A$ be any annulus component of $W' \cap \Theta$. By Proposition~\ref{spanningdisk}, $A$ has a spanning disk $D$ which we fix in $M$. If $\mathcal{B}$ is a 3-ball containing $\tau$ and intersecting $S'$ in a single circle, then by isotoping $\mathcal{B}$ we can move $\tau$ so that the core of $\tau$ is a properly embedded arc in $D$. Then using $D$, $\tau$ can be isotoped to be adjacent to $A$. If $\tau$ is in $W'$, the same argument holds not for $\tau$ but for a dual tube $\tau'$, where the disk in $V'$ giving the triviality of $\tau$ becomes the cocore of $\tau'$.

\begin{figure}[h] 
   \centering
   \includegraphics[width=3in]{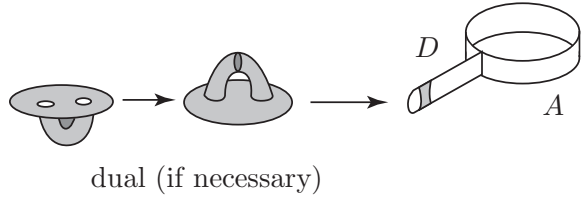} 
   \caption{A trivial tube is adjacent to any component of $S'$ cut along $\Theta'$.}
   \label{fig:tubemigration}
\end{figure}

\end{proof}

Now we show that stabilizing \hs \ in a horizontal component yields tubes which can be isotoped into other components of $M$ cut along $\Theta'$. Let $(X, \partial X, \emptyset)$ be a partition of $X$. Set $a(X) = a(X, \partial X, \emptyset) = g_h (X, \partial X, \emptyset) - 1$. Note that if $G$ is a fundamental 2-complex of $X$ then $|G| = a(X)$ (the Heegaard genus of $(X, \partial X, \emptyset)$ is $2g + m + k - 1$ if $k \geq 1$ and $2g + m$ if $k =0$; see Remark~\ref{amalgamationgenusofseifertfiberedspaces}).  

\begin{lemma}
\label{tubemigrationlemma}
Let $X$ be a horizontal component of $M$ cut along $\Theta'$ and suppose that for some component $Z_i$ of $V \cap X$, there is a tube in $W' \cap X$ (resulting from either stabilization or isotopy) which is adjacent to a component of $Z_i \cap \partial X$. Then $S' \cap Z_i$ can be isotoped so that all but at most $a(X)$ tubes are isotoped into neighboring components of $M$ cut along $\Theta'$. Moreover, if $S'$ is isotoped this way in every component of $V \cap X$, then at most $a(X)$ tubes remain in $X$.
\end{lemma}

\begin{proof}

Suppose first that $V \cap X$ has only one component. Let $S' \cap X$ be obtained from $S \cap X$ by adding a tube to $W$ that is adjacent to a component of $V \cap \partial X$. Let $A_1$ denote a component of $W' \cap \partial X$ which has a spanning disk that meets $X$ only in a spanning arc of $A_1$, and let $T_1$ be the component of $\partial X$ containing $A_1$. Now let $G$ be a fundamental 2-complex for $X$ built at $T_1$. Isotope $S' \cap X$ to be in tube position with respect to $G$. 

As the number of components of $G$ is $a(X)$, isotoping $S' \cap X$ to be in tube position gives at least $a(X)$ tubes in $V'$ with cores lying on components of $G$, and any two tubes on the same component of $G$ can be taken to be adjacent. Moreover, if there are at least $a(X) + 1$ tubes, then some component of $G$ has two adjacent tubes, and they can be chosen so that the annulus $A_1$ is between the feet of the tubes on one end, as in Figure~\ref{fig:tubeslide}. Let $A_2$ denote the component of $W' \cap \partial X$ between the feet on the other ends of the tubes, and let $T_2$ be the component of $\partial X$ containing $A_2$ (possibly $T_1 = T_2$). Performing a tube slide, we may isotope one of the tubes $\tau$ to be adjacent to $A_2$ in the component $X'$ of $M$ cut along $\Theta'$ across $T_2$ from $X$. 

If $A_1 = A_2$, then $\tau$ is trivial due to the spanning disk for $A_1$ used in the tube slide. Thus, applying Lemma~\ref{trivialtubelemma}, $\tau$ (or a dual tube) may be isotoped to be adjacent to another component which we call $A_2$ of $W' \cap \partial X$. If $A_1 \neq A_2$, then $\tau$ forms a spanning disk for $A_2$ in $X'$. In either case, there is a spanning disk formed for some component $A_2$ of $W' \cap \partial X$ that is not $A_1$. Note that if there is no other such component, then the tube $\tau$ is trivial. 

If there are at least $a(X)+1$ remaining tubes in $X$, then the spanning disk obtained for $A_2$ is one that can be used for a new tube slide (modifying the choice of the tubes obtained from isotoping $S' \cap X$ into tube position if necessary). Thus, by repeating the above process, all but $a(X)$ tubes (one for each component of $G$) can be isotoped out of $X$ by tube sliding. 
\smallskip

Now suppose that $V \cap X$ has $n$ components where $n > 1$, and let $Z_1$ be a component of $V \cap X$. Add a tube to $W \cap X$ (either by stabilization or isotopy) so that it is adjacent to $Z_1 \cap \partial X$. Suppose that no other component of $V \cap X$ has such a tube. The resulting surface $S' \cap X$ consists of a connected component $S' \cap Z_1$ along with $2(n-1)$ remaining connected horizontal essential surfaces $\hat{S'}$. Let $W_1$ be a component of $W \cap X$ such that $Z_1 \cap W_1$ is a component of $S \cap X$. Then some component $A_1$ of $W_1 \cap \partial X$ has a spanning disk $D_1$ meeting $W_1$ only in a spanning arc of $A_1$. As before, let $T_1$ denote the component of $\partial X$ containing $A_1$ and let $G$ be a fundamental 2-complex of $X$ built at $T_1$. Isotope $S' \cap Z_1$ to be in tube position with respect to $G$. 

Let $\tau$ be a tube whose foot is isotopic into $A_1$ (such a tube exists, since $S' \cap Z_i$ in tube position implies that for every component $A$ of $Z_1 \cap \partial X$ there is some tube obtained from ambient 1-surgery along an arc in $\mathcal{T}$ with at least one endpoint on an annulus parallel in $X$ to $A$). Using $D_1$, isotope $\tau$ via a tube slide to be adjacent to a spanning arc of the annulus $A_2$, as in Figure~\ref{fig:tubeslide2}. Thus $\tau$ is adjacent to $A_2$ in some component $X'$ of $M$ cut along $\Theta'$, and hence produces a spanning disk for $A_2$ in a neighboring component of $M$ cut along $\Theta'$. 

Note that it could be the case that $D_1$ runs along $\tau$ not in $W_1$ but on the other side of $Z_1$ in a component $W_1'$ of $W \cap X$ that also meets $Z_1$ nontrivially. If this is the case, then some component $A_1'$ of $W_1' \cap \partial X$ has a spanning disk given by a component of $D_1 - W_1'$ which does not run along $\tau$ and meets $W_1'$ only in a spanning arc of $A_1'$ (see Figure~\ref{fig:T2timesIspanningdisk} for a schematic of this phenomenon). If this is the case, take $W_1'$ as the component $W_1$, $A_1'$ as $A_1$ and the spanning disk component of $D_1 - W_1'$ as $D_1$ and apply the above argument. 

As before, if $A_1 = A_2$ then $\tau$ is trivial via the spanning disk $D_1$ for $A_1$. Thus applying Lemma~\ref{trivialtubelemma}, there is a spanning disk available for a new component $A_2$ of $W_1 \cap X$. If $A_1 \neq A_2$, then $\tau$ gives a spanning disk of the original $A_2$. We repeat the process until every tube formed by $S' \cap Z_1$ is isotoped out of $X$ to be adjacent to a component of $W \cap \partial X$. The only remaining components of $S' \cap Z_1$ after this isotopy are annuli which are parallel to the annuli in $Z_1 \cap \partial X$. Note also that in forming the spanning disks used to slide tubes of $S' \cap Z_1$, each component of $W_1 \cap \partial X$ obtained a spanning disk in a neighboring component of $X$, either from an adjacent tube or a spanning disk of the annulus that existed before isotopy. Further isotopy of the tubes outside of $X$ will not effect the existence of these spanning disks (we will discuss this formally in the next section). 

Now suppose $Z_i$ is some component of $V \cap X$, and suppose that in the components $Z_1, \ldots, Z_{i-1}$, $S'$ has been isotoped as above. Assume that a tube is isotoped into $Z_i$ and $S' \cap Z_i$ is isotoped into tube position. Tube slides for tubes formed by $S' \cap Z_i$ can now be performed using spanning disks for components of $W_i \cap \partial X$ along with previously established spanning disks from the above procedure (see Figure~\ref{fig:twocomponentslide}). Repeating this process eventually reduces the proof to the case where $S \cap X$ has two components. Thus the result follows.  

\begin{figure}[h] 
   \centering
   \includegraphics[width=1in]{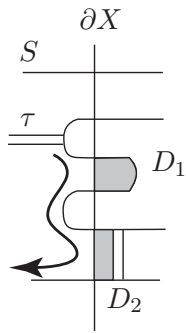} 
   \caption{A tube slide using two spanning disks. The disk $D_2$ is formed from tube sliding a previous component of $S' \cap X$.}
   \label{fig:twocomponentslide}
\end{figure}

\end{proof}

\section{Spanning disks and tube isotopies}
In the previous section we showed that if a tube is added to a component of $V' \cap X$ where $X$ is a horizontal component, then all but at most $a(X)$ tubes can be isotoped out of $X$ into neighboring components. The remainder of our discussion centers on how to isotope those tubes once they are out of $X$. This first requires a discussion of spanning disks of components of $\Theta'$ cut along $S'$. 

\subsection{Spanning disks}
\label{sec:section7}

Suppose that $A$ is a component of $W \cap \Theta'$ (before stabilizing). By Proposition~\ref{spanningdisk}, $A$ has a spanning disk $D$. After applying an innermost disk, outermost arc argument, $D$ can be assumed to intersect components of $M$ cut along $\Theta'$ in the following ways:
\medskip

\begin{itemize}
\item If $X$ is a horizontal component and $Z = F \times I$ is a component of $W \cap X$ where $F$ is a punctured surface, then $D \cap Z$ is a disjoint union of disks of the form $\delta_i \times I$, where the arcs $\delta_i$ are disjoint essential arcs on $F$. 
\medskip

\item If $Y$ is a vertical component, and $U$ is a component of $W \cap Y$ that does not contain an exceptional fiber or a component of $\partial M$ (so $U$ is a solid torus), then $D \cap U$ consists of a disjoint union of meridian disks of $U$. (Note that if $U$ contains an exceptional fiber or component of $\partial M$, then $D$ is assumed not to intersect $U$.)
\medskip

\item If $Q$ is a $T^2 \times I$ component, then each component of $W \cap Q$ is a solid torus and hence $D \cap Q$ is a disjoint union of meridian disks.
\medskip

\item If $N$ is the active component, then $D \cap W \cap N$ is a disjoint union of boundary compressing disks for $S \cap N$ in $W \cap N$, and possibly meridian disks of any solid tori components of $W \cap N$ if $N$ is an edge manifold and is aligned (see Figure~\ref{fig:option1}).
\end{itemize} 
\bigskip

\noindent We shall say that a spanning disk isotoped in this way is in {\em standard position}. 
\medskip

Let $X$ be a horizontal component, and let $Z$ be a component of $V \cap X$. Consider $Z = F \times I$ where $F$ is a punctured surface. Note that spanning disks \underline{in $W$} for components of $\Theta' \cap W$ will intersect $Z$ in a (possibly empty) disjoint union of essential arcs in $F \times \{0\}$ or $F \times \{1\}$. In particular, if \hs \ is stabilized by adding a tube $\tau$ to $W$ so that $\tau$ is adjacent to some component of $Z \cap \Theta'$, then the cocore of $\tau$ can be assumed to be disjoint from the arcs on $F \times \{0\}$ and $F \times \{1\}$. In other words, spanning disks for $W$ persist under stabilization of \hs \ in $Z$.

\begin{figure}[h] 
   \centering
   \includegraphics[width= 4 in]{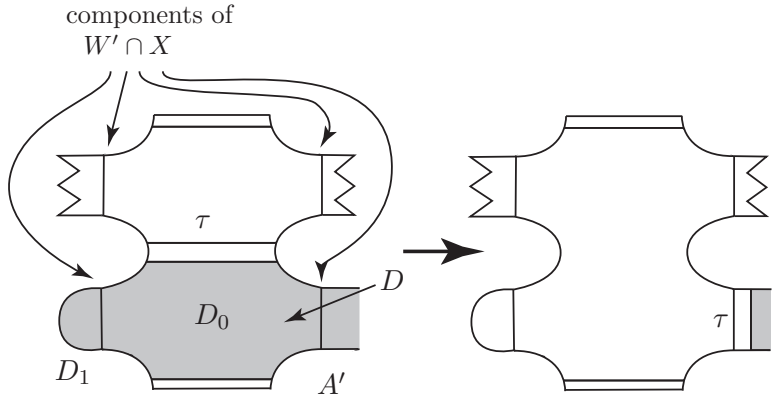} 
   \caption{A tube slide cuts off the spanning disk $D$. The disk $D_1$ is the spanning disk used for the tube slide as in Definition~\ref{tubeslidedefinition}.}
   \label{fig:flattubeslide}
\end{figure}

Next we observe what happens to a spanning disk via tube sliding. Suppose that $A$ is a component of $W \cap \Theta'$ and $D$ is a spanning disk for $A$ in standard position. Suppose that $D \cap Z \neq \emptyset$. Then after stabilizing \hs \ in $Z$ and isotoping $S' \cap Z$ into tube position, $D \cap S' \cap Z$ is a disjoint union of essential arcs, each arc running over some tube once. In particular, each component of $X \cap D$ can be taken to be a component of $G$ cut along $S'$ (or a parallel copy of such a component), where $G$ is a fundamental 2-complex for $X$. Suppose that $\tau$ and $\tau'$ are tubes meeting the boundary of some component $D_0$ of $D \cap X$ nontrivially. Perform a tube slide of $\tau$ so that $\tau$ becomes adjacent to a component $A'$ of $D_0 \cap \partial X$. Thus we may replace $D$ with the component of $D$ cut along $A'$ not containing $D_0$ (or component of $D$ not containing $D_0$ and parallel copies of $D_0$ if $D$ ran between $\tau$ and $\tau'$ more than once) as a spanning disk for $A'$. See Figure~\ref{fig:flattubeslide}.
\medskip

\subsection{Tube isotopies}

We have already described two ways of isotoping tubes, namely via isotoping trivial tubes as in Lemma~\ref{trivialtubelemma}, and via tube slides as in Lemma~\ref{tubemigrationlemma}. In particular, both of these types of isotopy can be assumed to isotope tubes along components of spanning disks. Now we add a third type of isotopy to our list.

\begin{figure}[h] 
   \centering
   \includegraphics[width=3.5in]{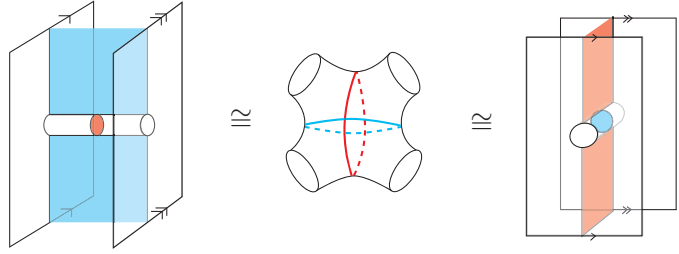} 
   \caption{A dual tube swap performed between two annuli in a solid torus connected by a tube.}
   \label{fig:dualtubeswap}
\end{figure}

\begin{definition}
\label{dualtubeswap}
\rm Suppose that $U$ is a solid torus in $M$ and that $S' \cap U$ is obtained from two properly embedded longitudinal annuli by ambient 1-surgery along an arc parallel into $\partial U$ which connects them. Let $\tau$ be the tube obtained from the ambient 1-surgery, and suppose $\tau$ is in $V'$. Then a {\em dual tube swap} of $\tau$ is an isotopy of $S' \cap U$ in $U$ yielding a tube $\tau'$ in $W'$ connecting two longitudinal annuli in $U$. See Figure~\ref{fig:dualtubeswap}.
\end{definition}

\begin{remark}
\label{typesofisotopyremark}
\rm These three types of isotopy -- tube slides, trivial tube isotopies, and dual tube swaps -- in addition to isotoping tubes guided by spanning disks, will be the ways that we isotope tubes once they are out of horizontal components. Note that if $\tau$ is isotoped to be adjacent to some component $A$ of $W' \cap \Theta'$, say, then if $\tau$ is further isotoped via spanning disks, trivial tube isotopies and/or dual tube swaps, so long as $\tau$ does not go through $A$ there will continue to exist a spanning disk for $A$. That is: $(1)$ if $\tau$ is isotoped via spanning disks or tube slides then a spanning disk for $A$ will continue to run along $\tau$ once, $(2)$ if $\tau$ is trivial then the disk giving the triviality of $\tau$ is a spanning disk for $A$ disjoint from $\tau$, and $(3)$ if a dual tube swap of $\tau$ is performed then there is a spanning disk for $A$ which is disjoint from the dual tube formed. Note also that the components $\Theta'$ cut along $S'$ can be assumed to be unchanged after each type of isotopy. 
\end{remark}

\section{Mutual Separability, $T^2 \times I$ components and edge manifold active components}
\label{sec:section8}

We have seen in Section~\ref{sec:horizontalcomponents} that stabilizing \hs \ in a horizontal component yields tubes which can be isotoped into neighboring components. Using the types of isotopy described in the previous section, those tubes can further be isotoped in the neighboring components, which can be $T^2 \times I$ components, the active component, or vertical components. In the following sections we shall deal with each scenario in turn. In this section we discuss what happens when tubes are isotoped into a $T^2 \times I$ component $Q$ or into the active component $N$ if $N$ is an edge manifold (as mentioned in Section~\ref{sec:horizontalcomponents}, if $N$ is a Seifert fibered component with $S \cap N$ psuedohorizontal, then we treat it as a horizontal component after isotopy -- see Remark~\ref{pseudohorizontalremark}). Note that in both situations, the component $Q$ or $N$ is homeomorphic to $T^2 \times I$.

\subsection{Mutual separability and blocked annuli} 

Before discussing the isotopies of tubes in components homeomorphic to $T^2 \times I$, we need to consider the general implications of $\Theta'$ being mutually separating.

\begin{definition}
\label{blockeddefinition}
\rm Let $X$ be a component of $M$ cut along $\Theta'$. Suppose that a Heegaard splitting \hhs \ of $(M, \partial_1 M, \partial_2 M)$ is isotoped so that $S'$ meets each component of $\Theta'$ in simple closed curves essential on both surfaces. A component $A$ of $\partial X$ cut along $\Theta'$ is called {\em blocked in $X$} if there is a spanning disk $D$ of $A$ contained entirely in $X$. Otherwise $A$ is {\em unblocked in $X$}.
\end{definition}

\begin{figure}[h] 
   \centering
   \includegraphics[width=3in]{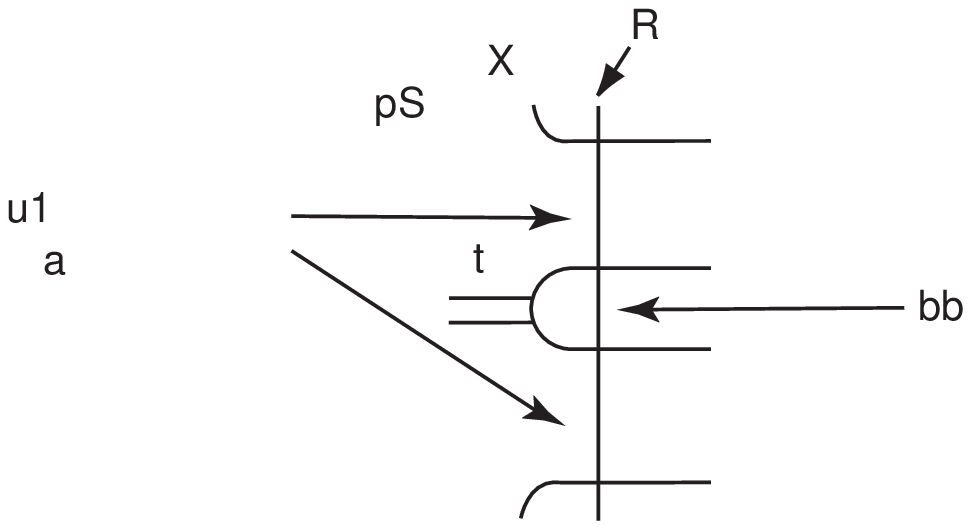} 
   \caption{Blocked and unblocked annuli in $\partial X$ cut along $S'$.}
   \label{fig:blockedandunblocked}
\end{figure}

Recall that by construction, $\Theta'$ is mutually separating. Thus, the components of $M$ cut along $\Theta'$ can be partitioned into two manifolds $M_1$ and $M_2$ such that every component of $\Theta'$ has a neighborhood that meets $M_1$ and $M_2$ nontrivially. To establish that \hs \ stabilizes to an amalgamation along $\Theta'$, Lemma~\ref{theamalgamationlemma} requires that every component of $V' \cap \Theta'$ is blocked in $M_1$ and every component of $W' \cap \Theta'$ is blocked in $M_2$, or vice versa.  Note that if $X$ is a horizontal component and every component of $V' \cap X$ is isotoped into tube position as per Lemma~\ref{tubemigrationlemma}, then it is easy to see that every component of $V' \cap \partial X$ is blocked in $X$. 

\begin{remark}
\label{startingstabilizationremark}
\rm Suppose that $X$ is a horizontal component of $M$ cut along $\Theta'$ such that \hs \ is stabilized in $V \cap X$. Then this determines that the blocked components of $\partial X$ cut along $\Theta'$ are in $V' \cap \partial X$. If $X \subset M_1$, this implies that every component $X'$ of $M$ cut along $\Theta'$ must have blocked annuli in $V' \cap \partial X'$ if $X' \subset M_1$ and in $W' \cap \partial X'$ if $X' \subset M_2$ in order to apply Lemma~\ref{theamalgamationlemma}. In particular, this determines a set of ideal partitions for components of $M$ cut along $\Theta'$, namely the components of $X' \cap \Theta'$ must all be contained in $\partial_1 X'$ when $X' \subset M_1$, and in $\partial_2 X'$ when $X' \subset M_2$. Thus we will choose to stabilize \hs \ in either $V \cap X$ or $W \cap X$, depending on which gives rise to the ideal set of partitions of $M$. 
\end{remark}

\subsection{$T^2 \times I$ components}

Let $Q$ be a $T^2 \times I$ component of $M$ cut along $\Theta'$ which is not the active component. Then $S \cap Q$ consists of essential annuli, each of which has its boundary components on different components of $\partial Q$.  Note that $Q$ cannot meet vertical components at both its boundary components, for this would imply that the fiberings of the two vertical components and $Q$ line up so that $\partial Q$ would not have been chosen as components of $\Theta'$. Let $X_1$ and $X_2$ be the components of $M$ cut along $\Theta'$ meeting $\partial Q$ (possibly $X_1 = X_2$). We show that a tube can often be passed through $Q$ and moreover, if $S'$ can be isotoped so that components of $V' \cap \partial X_i$ are blocked in $X_i$, $i = 1,2$, then then the components of $W' \cap \partial Q$ can be blocked in $Q$.

\begin{lemma} 
\label{T2timesIlemma}
Let $\tau \subset V'$ be a tube isotoped into a $T^2 \times I$ component $Q$. Then if $S' \cap Q$ has essential annular components disjoint from $\tau$, after a dual tube swap of $\tau$ the dual tube $\tau'$ can be isotoped into $X$ or $X'$ so that $\tau'$ is trivial or adjacent to a component of $V' \cap \partial Q$. Moreover, if every component of $V' \cap \partial Q \cap \partial X_i$ is blocked in $X_i$, $i = 1, 2$, then after isotoping a tube in $V'$ into $Q$, $S' \cap Q$ can be isotoped so that every component of $W' \cap \partial Q$ is blocked in $Q$. 
\end{lemma}

\begin{proof}
If $\tau \subset V'$ is isotoped into $Q$, then a dual tube swap yields a dual tube $\tau' \subset W'$. Assume that there exist components of $S' \cap Q$ which are essential annuli (thus disjoint from $\tau'$). Then the feet of $\tau'$ are isotopic into one of four possible components $A_1, A_2, A_3$ or $A_4$ of $V' \cap \partial Q$. We note that there has to be some $A_i$ with a spanning disk disjoint from $\tau'$ (if $\tau$ was tube slid out of $X_i$ where $X_i$ is a horizontal component, this is obvious). To see this, take a spanning disk $D_1$ of $A_1$, say. Such a spanning disks exists by Lemma~\ref{spanningdisk}. If the disk $D_1$ is not disjoint from $\tau'$, then some outermost disk component $D$ of $D_1$ cut along $\cup_{i=2}^4 A_i$ is. This disk $D$ is our desired spanning disk (see Figure~\ref{fig:T2timesIspanningdisk}). 

  \begin{figure}[h] 
   \centering
   \includegraphics[width= 2.4in]{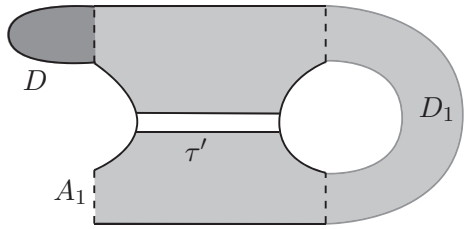} 
   \caption{A spanning disk for a component of $V' \cap \partial Q$ that is disjoint from $\tau'$.}
   \label{fig:T2timesIspanningdisk}
\end{figure}

Thus $\tau'$ can be isotoped out of $Q$ by an isotopy guided by $D$ so that $\tau'$ becomes adjacent to a component of $V' \cap \partial Q$. Note that this introduces a pair of components of $S' \cap Q$ that are parallel to components of $W' \cap \partial Q$ giving spanning disks in $Q$ for those components. They are next to components of $V' \cap \partial Q$ which have spanning disks, namely $D$ and the spanning disk in $X_2$ formed by $\tau'$ after isotopy. If $\tau'$ is trivial then then disk giving its triviality is the spanning disk for that component of $W' \cap \partial Q$. By Remark~\ref{typesofisotopyremark} the existence of these spanning disks will persist under additional isotopy of tubes.  

 \begin{figure}[h] 
   \centering
   \includegraphics[width= 2.75in]{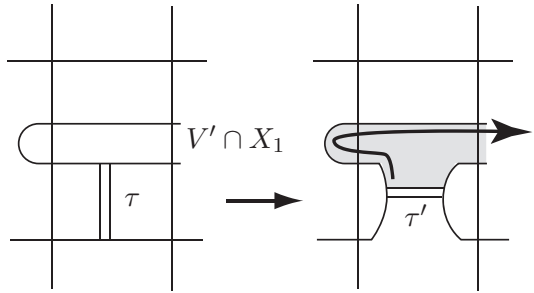} 
   \caption{Isotopy of $\tau'$ to be adjacent to a component of $V' \cap \partial Q$ in $X_i$.}
   \label{fig:T2_swap2}
\end{figure}

The second assertion of the lemma follows from the fact that after isotoping $\tau$ into $Q$ and then dual tube swapping, the existence of spanning disks for each component of $V' \cap \partial Q$ outside of $Q$ implies that the dual tube $\tau'$ is trivial. 

\begin{figure}[h] 
   \centering
   \includegraphics[width=4in]{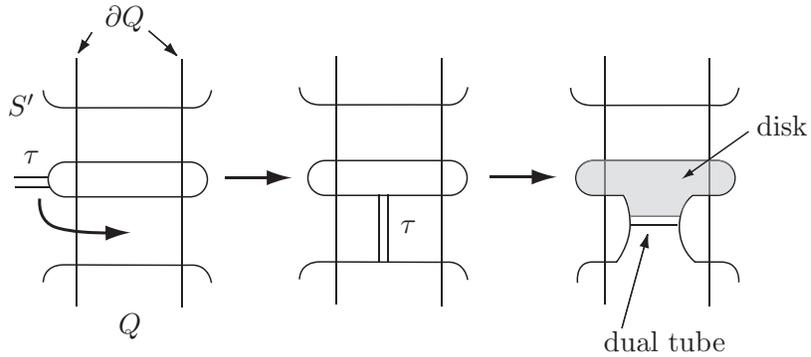} 
   \caption{Dual tube swapping $\tau$ in $Q$ yields a trivial tube.}
   \label{fig:trivialinQ}
   \end{figure}

By Lemma~\ref{trivialtubelemma}, $\tau'$ can thus be isotoped into $Q$ to connect other essential annular components of $S \cap Q$. After repeating this argument until all essential annular components of $S' \cap Q$ have been eliminated, $S' \cap Q$ is such that every component of $W' \cap \partial Q$ is blocked in $Q$ (see Figure~\ref{fig:trivialinQ}). 

\begin{figure}[h] 
   \centering
   \includegraphics[width=5in]{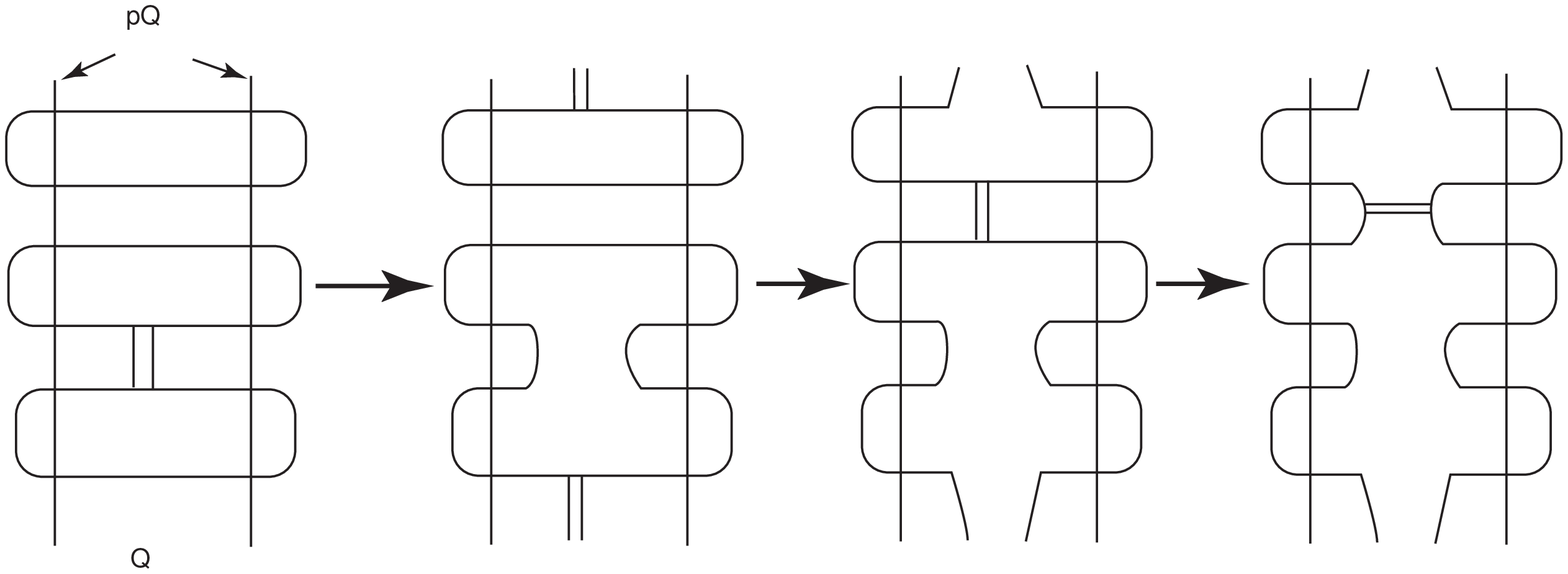} 
   \caption{Isotoping $S' \cap Q$ so that every component of $W' \cap \partial Q$ is blocked in $Q$.}
   \label{fig:trivialinQ2}
\end{figure}

\end{proof}

Note that if $S' \cap Q$ is isotoped as above, then any additional tubes isotoped into $Q$ can be isotoped through $Q$ or are trivial. Also note in order to isotope $S' \cap Q$ so that every component of $W' \cap Q$ is blocked in $Q$, $Q$ requires 1 tube isotoped into it.

\subsection{The active component as an edge manifold}
Suppose the active component $N$ is an edge manifold. Then either $N$ is a $T^2 \times I$ component or is part of a Seifert fibered component of $M$ cut along $\Theta'$. In each of these situations we obtain an analogue of Lemma~\ref{T2timesIlemma} for $N$. 

\begin{lemma}
\label{activecomponentlemma}
Suppose the active component $N$ is an edge manifold.
\begin{itemize}
\item[1.] If $N$ is a $T^2 \times I$ component and $\tau$ is a tube in $V'$ isotoped into $N$, then either $\tau$ is trivial or after a dual tube swap it can be isotoped to be adjacent to a component of $V' \cap \partial N$ in a neighboring component of $N$. Moreover, if every component of $V' \cap \partial N$ is blocked in neighboring components of $N$, $S' \cap N$ can be isotoped so that every component of $W' \cap \partial N$ is blocked in $N$. 
\medskip

\item[2.] If $N$ is contained in a Seifert fibered component $X$ ({\em i.e.}~ $N$ intersects $\Theta'$ in only one component $T$), then a tube $\tau$ in $V'$ isotoped into $N$ is either trivial or can be isotoped through $N$ to be adjacent to a component of $W' \cap \partial N$. Moreover, if every component of $V' \cap (\partial N - T)$ is blocked in $\overline{X-N}$, then $S' \cap N$ can be isotoped so that some component of $W' \cap T$ is blocked in the component of $M$ cut along $\Theta'$ meeting $N$ at $T$.
\end{itemize}
\end{lemma}
\begin{proof} 
Suppose first that $N$ is a $T^2 \times I$ component. If $S' \cap N$ is a toggle, then any tube isotoped into $N$ is trivial. Moreover, if the components of $V' \cap \partial N$ are blocked in neighboring components of $N$, then it is easy to see that the components of $W' \cap \partial N$ are blocked in $N$.

Suppose that $S' \cap N$ is aligned. If $\tau$ is isotoped into $N$ such that $\tau$ connects two essential annular components of $S' \cap N$, then as in the proof of Lemma~\ref{T2timesIlemma}, after performing a dual tube swap the resulting dual tube $\tau'$ can be isotoped to be adjacent to a component of $V' \cap \partial N$ in a neighboring component of $N$. Also as in the proof of Lemma~\ref{T2timesIlemma}, this results in components of $S' \cap N$ which are parallel to components of $W' \cap \partial N$, with components of $V' \cap \partial N$ next to them having spanning disks outside of $N$. If another tube is isotoped into $N$ with one foot on such a boundary parallel component, then the aforementioned disks either give the tube as trivial or provide an isotopy of the foot of the tube onto a component of $S' \cap N$ which is not boundary parallel. Thus the above argument may be repeated.

If at least one of the feet $\tau$ lies on the 4-times punctured sphere component of $S' \cap N$, then either $\tau$ is trivial or performing a dual tube swap allows the dual tube $\tau'$ to be isotoped either to be trivial or to be adjacent to a component of $V' \cap \partial N$ in a neighboring component of $N$. This uses the fact, similar to what was shown in the proof of Lemma~\ref{T2timesIlemma}, that some spanning disk for a component of $V' \cap \partial N$ into which a foot of $\tau'$ is isotopic is disjoint from $\tau'$. An example of such a situation is given in Figure~\ref{fig:tubethroughactivecomponent}.

\begin{figure}[h]
   \centering
   \includegraphics[width=3in]{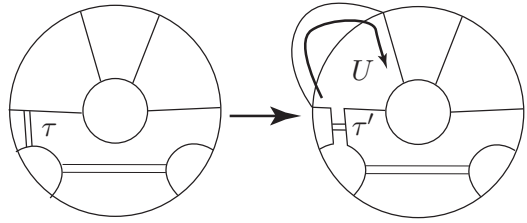} 
   \caption{A dual tube slide yields the tube $\tau'$. As in the proof of Lemma~\ref{T2timesIlemma} some component of $U \cap \Theta'$ has a spanning disk disjoint from $\tau'$ that can be used to isotope $\tau'$ through $N$.}
   \label{fig:tubethroughactivecomponent}
\end{figure}

If every component of $V' \cap \partial N$ is blocked in neighboring components of $N$, then it is straightforward to show (possibly using dual disk swaps and trivial tube isotopies) that $S' \cap N$ is isotopic to a disjoint union of boundary parallel annuli with a tube connecting a pair of such annuli with boundary components on different components of $\partial N$ (as in Figure~\ref{fig:trivialinQ2}). Thus the first conclusion of the lemma follows.
\medskip

Now suppose that $N$ is contained in some Seifert fibered component $X$. Let $X_0 = \overline{X - N}$ and set $T' = N \cap X_0$ and $T = \partial N - T'$. As before, if $S' \cap N$ is a toggle then every tube isotopic into $N$ is trivial. If $S' \cap N$ is aligned then it is easy to see that any tube in $V'$, say, isotoped into $N$ is trivial or can be isotoped through $N$ (possibly using a spanning disk of an annulus as in Figure~\ref{fig:tubethroughactivecomponent}) to be adjacent to a component of $W' \cap N$ in a neighboring component of $N$. 

Suppose that every component of $V' \cap T$ is blocked in $X_0$. If $S' \cap N$ is a toggle, then $S' \cap X$ can be isotoped as in the proof of Lemma~\ref{theamalgamationlemma} so that $S'$ intersects $T'$ in $T'$ with open disk removed. In particular, this allows for a dual tube swap of the tube in the toggle, which can be isotoped through $T$ to be adjacent to a component of $W' \cap T$ in the component of $M$ cut along $\Theta'$ that meets $N$ at $T$ (see Figure~\ref{fig:toggletubeswap}).

\begin{figure}[h]
   \centering
   \includegraphics[width=3in]{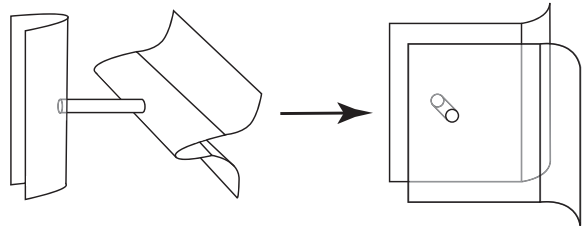} 
   \caption{Isotoping a toggle to yield a tube adjacent to a component of $\Theta'$ cut along $S'$.}
   \label{fig:toggletubeswap}
\end{figure}

If $S' \cap N$ is aligned, then either a tube can be isotoped directly across a component of $W' \cap T$ (as in Figure~\ref{fig:activetubecontribution}), or if such a tube is isotopic into $T$, then (possibly after a dual tube swap) there is a trivial tube which can be isotoped via Lemma~\ref{trivialtubelemma} to be adjacent to a component of $W' \cap T$. 

\begin{figure}[h]
   \centering
   \includegraphics[width=3in]{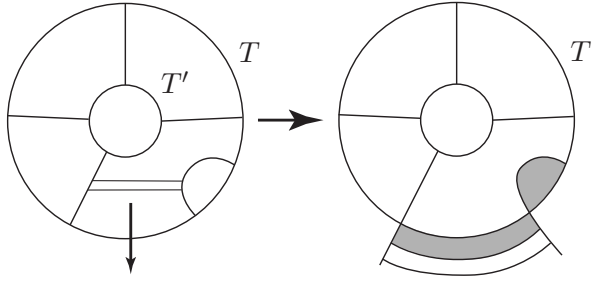} 
   \caption{Isotoping $S' \cap N$ to give a spanning disk on the other side of $T$ in the case that $S' \cap N$ is aligned.}
   \label{fig:activetubecontribution}
\end{figure}

\end{proof}

The importance of Lemma~\ref{activecomponentlemma} is twofold. First, it shows that if $N$ is an edge manifold, then tubes can be appropriately passed through $N$ and into other components of $M$ cut along $\Theta'$. Secondly, it shows that if $S'$ is isotoped so that components of $V' \cap \partial N$, say, are blocked in neighboring components to $N$, then one less tube is needed from tube isotopies to block all the components of $W' \cap \partial N$. If $N$ is a $T^2 \times I$ component, then the active component supplies the tube in $N$ that otherwise would need to be isotoped into a $T^2 \times I$ component to obtain spanning disks for every component of $W' \cap \partial N$, as in Lemma~\ref{T2timesIlemma}. If $N$ is part of a Seifert fibered component, then isotoping $S' \cap N$ as in Lemma~\ref{activecomponentlemma} shows that one less tube is needed for the components of $W' \cap \partial N \cap \Theta'$ to be blocked in the neighboring component to $N$. We will use this fact in Section~\ref{sec:proofofthemainresults}. 

\section{Vertical components}
\label{sec:section9}

Now we consider the case that a component of $M$ cut along $\Theta'$ is vertical. As we saw in the case that a tube was isotoped into a $T^2 \times I$ component, the tube could be passed on to a neighboring component of $Q$ unless it connected the last two essential annular components of $S' \cap Q$, in which case it stayed in $Q$. In this section we show a similar result for vertical components, in that a tube isotoped into a vertical component $Y$ can either be passed on to another component of $M$ cut along $\Theta'$ or remains in $Y$.

Let $Y$ be a vertical component of $M$ cut along $\Theta'$, and suppose that in every component $X$ of $M$ cut along $\Theta'$ that meets $Y$ (no such component can be vertical as this would imply that the fiberings of $Y$ and $X$ line up) $S' \cap X$ is isotoped so that every component of $V' \cap (\partial X \cap \partial Y)$ is blocked in $X$ ({\em e.g.}~ $X$ is horizontal and each component of $S' \cap X$ is in tube position). Set $\partial_1 Y = \partial_1 M \cap \partial Y$ and $\partial_2 Y = (\partial_2 M \cap \partial Y) \cup (\partial Y \cap \Theta).$ This gives a partition $(Y, \partial_1 Y, \partial_2 Y)$ of $Y$. We will call this partition {\em weighted} if $\partial_1 Y = \emptyset$ and $\partial_2 M \cap \partial Y \neq \emptyset$, and {\em unweighted} otherwise. For simplicity of notation, set $g_h (Y) = g_h (Y, \partial_1 Y, \partial_2 Y)$ and 
$a(Y) = g_h(Y) - 1$. Let $m_E$ be the number of components of $\partial M \cap Y$.

\begin{remark}
\label{amalgamationgenusofseifertfiberedspaces}
\rm If $Y$ is a Seifert fibered component of a graph manifold $M$ such that $Y$ has a partition $(Y, \partial_1 Y, \partial_2 Y)$ given as above, then 

$$g_h (Y) = \left\{ \begin{array} {ll}
	2g + m + k -1 & \mbox{if $k \geq 1$, or $m_E \geq 1$ and $(Y, \partial_1 Y, \partial_2 Y)$ is not weighted} \\
	2g + m & \mbox{if $k=0$, and $m_E = 0$ or $(Y, \partial_1 Y, \partial_2 Y)$ is weighted} 
	\end{array}	\right. $$

\noindent where $Y$ is given by the Seifert data $\left(g;m; \frac{\beta_1}{\alpha_1}, \ldots, \frac{\beta_k}{\alpha_k} \right).$

\end{remark}

If $k \geq 1$, or $m_E \geq 1$ and $(Y, \partial_1 Y, \partial_2 Y)$ is not weighted, then $g_h(Y)$ is the standard Heegaard genus for $Y$, as shown for example in \cite{Schultens4}. If $k=0$, and $m_E = 0$ or $(Y, \partial_1 Y, \partial_2 Y)$ is weighted, then all the boundary components of $Y$ lie in one of the compression bodies of the Heegaard splitting, hence $g_h(Y)$ is the standard Heegaard genus of $Y$ plus one. 

\begin{lemma}
\label{verticalimportlemma}
Let $Y$ be a vertical component of $M$ cut along $\Theta'$, and suppose that in every component $X$ of $M$ cut along $\Theta'$ which meets $Y$ (other than $Y$), $S' \cap X$ is isotoped so that every component of $V' \cap (\partial Y \cap \partial X)$ has a spanning disk in $X$.  Then after isotoping $a(Y)$ or $a(Y) + 1$ tubes into $Y$ every annulus component of $W' \cap (\partial Y - \partial M)$ has a spanning disk in $Y$. Moreover, in the case that $a(Y) + 1$ tubes are isotoped into $Y$, after isotopy the above conclusion holds and there is an additional tube in $Y$ which is trivial. 
\end{lemma}

\begin{proof} Suppose $Y$ has Seifert data

$$\left(g;m; \frac{\beta_1}{\alpha_1}, \ldots, \frac{\beta_k}{\alpha_k} \right).$$

Let $U$ be a component of $Y$ cut along $S$ (before stabilizing). Then $U$ is homeomorphic to $D^2 \times S^1$ or $T^2 \times I$. Let $\varepsilon ( U )$ be the number of exceptional fibers or the number of components of $\partial Y \cap \partial M$ contained in $U$. Note that $U$ cannot contain both an exceptional fiber of $Y$ and a component of $\partial Y \cap \partial M$, thus $\varepsilon (U) = 0$ or $1$. Let $u(U) = | U \cap (W' \cap (\partial Y - \partial M))|$ (note that $u(U) \geq 1$). Set $c(U) = u(U) + \varepsilon(U) - 1$.

Let $U_1, \ldots, U_r$ be the components of $W \cap Y$ and let $\Upsilon_1, \ldots, \Upsilon_s$ be the components of $V \cap Y$. Define

$$c(Y) = \sum_{i=1}^r c(U_i) + \sum_{i=1}^s c(\Upsilon_i).$$

Let $t$ denote the number of tubes needed to be isotoped into $Y$ such that every annulus component of $W' \cap (\partial Y - \partial M)$ has a spanning disk in $Y$. By assumption, every component of $V' \cap (\partial Y - \partial M)$ has a spanning disk in a neighboring Seifert fibered component to $Y$. 
\bigskip

\noindent {\bf Claim:} $t = c(Y)$ or $c(Y) + 1$. If $t = c(Y) + 1$, then either $(1)$ after sliding $t$ tubes into $Y$ and isotoping there is a trivial tube in $Y$, or $(2)$ $k=0$, and $m_E = 0$ or $(Y, \partial_1 Y, \partial_2 Y)$ is weighted. 
\bigskip

\begin{proof}
First note that if $c(U_i) = 0$ for some $i$, then $U_i \cap S$ is parallel to $ U_i  \cap \Theta'$. This contradicts the assumption that \hs \ is in active position, {\em i.e.}~all annuli in $\Theta' \cap W$ are essential in $W$. Thus $c(U_i) \geq 1$ for $1 \leq i \leq r$. Note that if $c(U_i)$ tubes are isotoped into $U_i$ such that each component of $U_i \cap (\partial Y - \partial M)$ has at most one adjacent tube, then each component of $W' \cap (\partial Y - \partial M)$ in $U_i$ has a spanning disk in $U_i$ (see Figure~\ref{fig:verticalslide}). Also note that if $c(U_i)$ tubes are isotoped into $U_i$, then for a given component of $U_i \cap (\partial Y - \partial M)$, possibly after isotopy, there is a tube adjacent to it. Moreover, if two tubes are isotoped to be adjacent to the same component of $W' \cap (\partial Y - \partial M)$, one of the tubes is trivial and can be isotoped via Lemma~\ref{trivialtubelemma} to be adjacent to another component of $W' \cap (\partial Y - \partial M)$ as needed. 
\medskip

\begin{figure}[h] 
   \centering
   \includegraphics[width=5in]{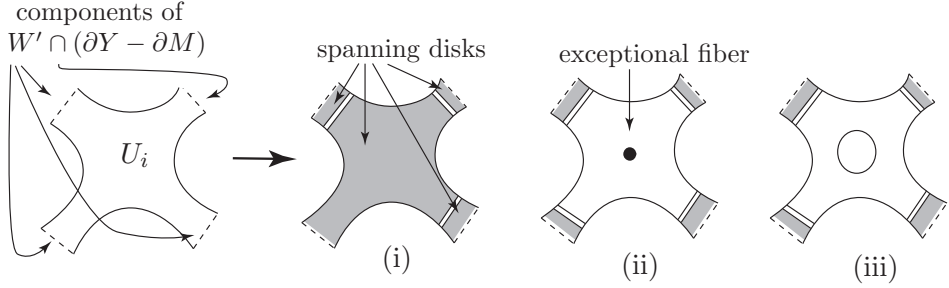} 
   \caption{The cases for $U_i$: (i) $\varepsilon(U_i) = 0$, (ii) $U_i$ contains an exceptional fiber, and (iii) $U_i$ contains a component of $\partial M$.}
   \label{fig:verticalslide}
\end{figure}

Now consider $\Upsilon_i$ for some $1 \leq i \leq s$. The quantity $c(\Upsilon_i) = 0$ or $-1$, depending on whether $\varepsilon(\Upsilon_i) = 1$ or $0$. If $c(\Upsilon_i) = 0$, then it contributes nothing to $c(Y)$, and therefore does not affect the calculation of $t$. If on the other hand $c(\Upsilon_i) = -1$, then we need to show that one of the $\sum_{i=1}^r c(U_i)$ potential tubes needed for each component of $W' \cap (\partial Y - \partial M)$ to have a spanning disk in $Y$ is in fact unnecessary, or if it is the case that all the components of $W' \cap (\partial Y - \partial M)$ already have spanning disks in $Y$, that a trivial tube is produced as in conclusion $(1)$ of the claim. 

So assume that $c(\Upsilon_i) = -1$. Thus $\Upsilon_i$ is a solid torus containing no exceptional fibers. As before, $|\Upsilon_i \cap \partial Y| \geq 2$ since $|\Upsilon_i \cap \partial Y| = 1$ implies that $V \cap \Theta'$ has an inessential annulus, contradicting the assumption that \hs \ is in active position. As $|\Upsilon_i \cap \partial Y| = |\Upsilon_i \cap S|$, this implies $\Upsilon_i$ meets components of $\cup_{i=1}^r U_i$ in at least 2 annuli.

Suppose that $U_l$ meets $\Upsilon_i$ in a component of $\Upsilon_i \cap S$. Assume that after isotoping tubes into $U_l$, there is a tube $\tau$ with a foot on that component of $\Upsilon_i \cap S$ (note that we can always assume this, by isotoping as many as $c(U_l)$ tubes into $U_l$). Suppose that the other foot of $\tau$ does not lie on $\Upsilon_i \cap S$, and this is the first such tube isotoped into $Y$ that satisfies this condition. Then by a dual tube swap, 
$\tau$ is replaced so that $\Upsilon_i$ is connected to some other $\Upsilon_j$ via two annuli (one of which is parallel into $\partial Y$) and a tube $\tau'$ connecting them. Since $\Upsilon_i$ contains no exceptional fibers or components of $\partial M$, a meridian disk of $\Upsilon_i$ gives that $\tau'$ is trivial (see Figure~\ref{fig:verticalslide2}). 

\begin{figure}[h] 
   \centering
   \includegraphics[width=4in]{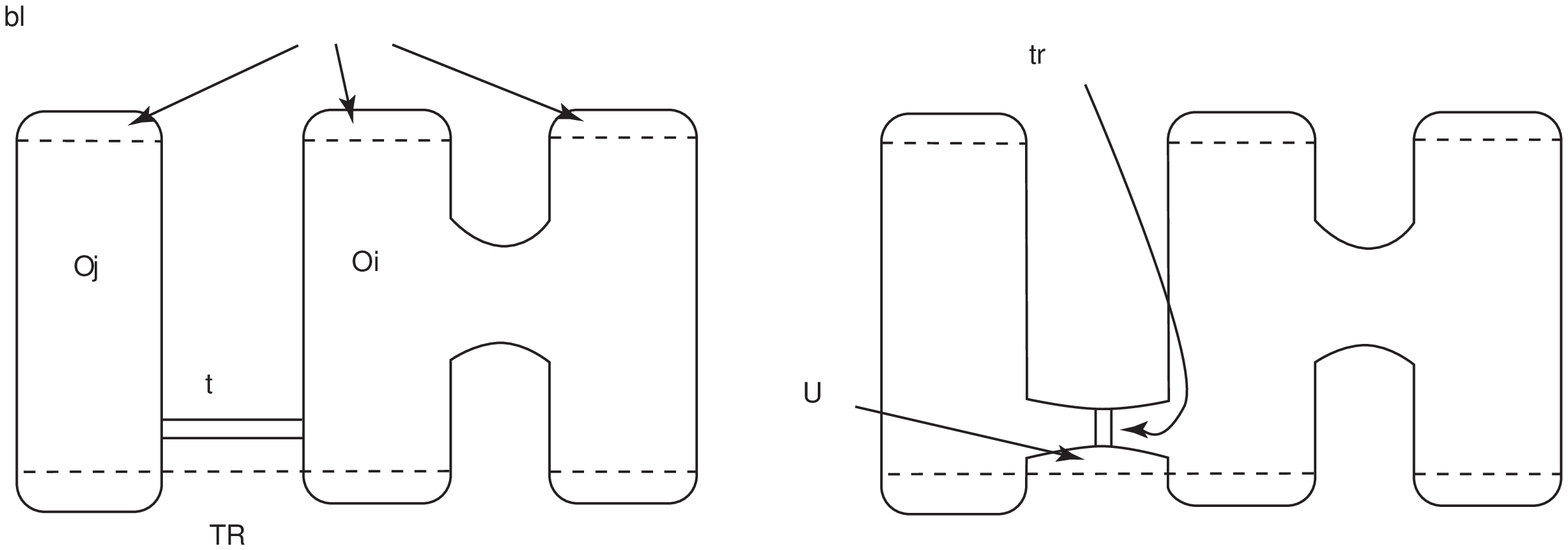} 
   \caption{A dual tube swap along $\tau$.}
   \label{fig:verticalslide2}
\end{figure}

As $\Upsilon_i$ is combined with $\Upsilon_j$, after redefining $U_1, \ldots, U_r$ and $\Upsilon_1, \ldots, \Upsilon_s$ the number $s$ goes down by one and the number $r$ goes up by one. Note that the new component $U$ of $W' \cap Y$ created in the dual tube swap is a solid torus with no exceptional fibers. A meridian disk of $U$ is thus a spanning disk for the corresponding component of $W' \cap \partial Y$, and moreover $c(U) = 0$ so it leaves the calculation of $c(Y)$ unaffected. 

By Lemma~\ref{trivialtubelemma}, isotope $\tau'$ to be adjacent to another component of $U_{l'} \cap (\partial Y - \partial M)$ (possibly in $l = l'$) as needed to form a spanning disk of that component in $Y$. Then the above argument can be repeated as necessary. Note that this lowers the number of tubes needed for components of $U_{l'} \cap (\partial Y - \partial M)$ to have a spanning disk in $Y$ by one, as desired. If, for every $1 \leq l \leq r$, every component of $U_l \cap (\partial Y - \partial M)$ already has a spanning disk in $Y$, then there is a trivial tube left over, and conclusion $(1)$ of the claim results.

The following argument shows that if $s > 1$, then for some $U_l$ we can always assume that after isotoping at most $c(U_l)$ tubes into $U_l$ there is a tube $\tau$ which for a fixed $i$ meets $\Upsilon_i$ at one of its feet and some $\Upsilon_j$, $j \neq i$, at the other. In particular, we show that if $s > 1$ then there is some component of $W' \cap (\partial Y - \partial M)$ which has one boundary component on $\Upsilon_i \cap S$ and the other on $\Upsilon_j \cap S$ for some $i \neq j$. To see this, note that $\cup_{i=1}^{r} \partial U_i - \partial M = (S \cap Y) \cup (W' \cap (\partial Y - \partial M))$. Assume that there is no component of $W' \cap (\partial Y - \partial M)$ intersecting in its boundary $\partial \Upsilon_i$ and $\partial \Upsilon_j$ for some $i \neq j$. Then we can set $S_i =  S \cap \partial \Upsilon_i$, $S_{\hat{i}} = \cup_{j \neq i} S \cap \partial \Upsilon_j$, and take $W'_{i}$ to be  the components of $W' \cap (\partial Y - \partial M)$ meeting $\Upsilon_i$ in both boundary components, and $W_{\hat{i}}$ to be the components of $W' \cap (\partial Y - \partial M)$ meeting $\Upsilon_j$ in both boundary components for $j \neq i$. Observe that $S \cap Y = S_i \sqcup S_{\hat{i}}$ and $W' \cap (\partial Y - \partial M) = W'_i \sqcup W'_{\hat{i}}$. Now, for any $U_l$, the component of $\partial U_l$ intersecting $\partial Y - \partial M$ is connected, hence $\partial U_l$ is contained in $S_i \cup W'_i$ or $S_{\hat{i}} \cup W'_{\hat{i}}$ for each $1 \leq l \leq r$. But then the sets $\Upsilon_i \cup \{ U_l | \partial U_l \subset S_i \cup W'_i \}$ and $\cup_{j \neq i} \Upsilon_j \cup \{ U_l | \partial U_l \subset S_{\hat{i}} \cup W'_{\hat{i}} \}$ disconnect $Y$, which is a contradiction. Thus there must be some $U_l$ such that after isotoping at most $c(U_l)$ tubes into $U_l$ there is a tube which is adjacent to a component of $\partial U_l \cap (\partial Y - \partial M)$ between $\Upsilon_i$ and $\Upsilon_j$. 
\smallskip

We are reduced to the case that after isotopy either $s \geq 1$ and $c(\Upsilon_i) = 0$ for all $1 \leq i \leq s$, or $s = 1$ and $c(\Upsilon_i) = -1$. In the former case, precisely $c(Y) = \sum_{i=1}^r c(U_i) + \sum_{i=1}^s c(\Upsilon_i) = \sum_{i=1}^r c(U_i)$ tubes are needed for every component of $W' \cap (\partial Y - \partial M)$ have a spanning disk in $Y$, and the conclusion of the claim is satisfied. Note that for some $c(\Upsilon_i)$ to equal zero, $\Upsilon_i$ must contain an exceptional fiber or a component of $\partial M$, implying that $k \geq 1$, or $m_E \geq 0$ and $(Y, \partial_1 Y, \partial_2 Y)$ is not weighted. 

For the final case, assume that $s = 1$ and $c(\Upsilon_i) = c(\Upsilon) = -1$. If $k \geq 1$, let $U_i$ be a component of $W \cap Y$ containing an exceptional fiber. Since $s=1$, $\partial U_i \cap \partial \Upsilon \neq \emptyset$, hence a tube isotoped into $U_i$ has both feet on $\partial \Upsilon$. 
Let $\tau$ be such a tube and assume first that both feet of $\tau$ lie on the same component of $S \cap \Upsilon$.

By a dual tube swap, we can take $\Upsilon$ to be a solid torus such that $S' \cap \Upsilon$ is obtained from a punctured torus in $\partial \Upsilon$ and the frontier of a neighborhood of the exceptional fiber in $U_i$ by ambient 1-surgery along an arc connecting the two. As Heegaard splittings of a solid torus are standard, this gives a trivial tube in $\Upsilon$ (see Figure~\ref{fig:redefineupsilon}). After applying this isotopy, the exceptional fiber initially in $U_i$ is now in $\Upsilon$ instead. Thus $c(U_i)$ has decreased by one and $c(\Upsilon)$ has increased by one, plus there is a trivial tube to be isotoped via Lemma~\ref{trivialtubelemma} to be adjacent to another component of $W' \cap (\partial Y - \partial M)$. Hence, the number of tubes needed to be isotoped into $Y$ to obtain a spanning disk in $Y$ of every component of $W' \cap (\partial Y - \partial M)$ is $\sum_{i=1}^r c(U_i) + c(\Upsilon)  + 1 = \sum_{i=1}^r c(U_i) = c(Y)$, or if there are already $c(Y)$ tubes in $Y$ so that every component of $W' \cap (\partial Y \cap \partial M)$ has a spanning disk in $Y$, there is a trivial tube as in conclusion $(1)$ of the claim. 

\begin{figure}[htbp] 
   \centering
   \includegraphics[width=3in]{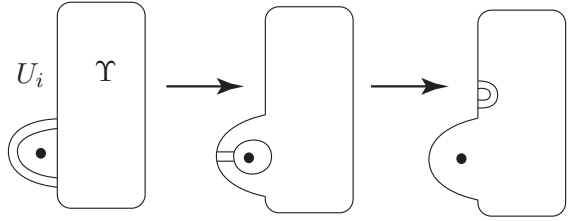} 
   \caption{Using a dual tube swap and isotopy to redefine $U_i$ and $\Upsilon$.}
   \label{fig:redefineupsilon}
\end{figure}

If the feet of $\tau$ lie on different components of $\Upsilon \cap S$, then we require that each component of $U_i \cap (\partial Y - \partial M)$ have a tube adjacent to it in $U_i$. Thus, as $c(\Upsilon) = -1$, we need $\sum_{i=1}^r c(U_i) + c(\Upsilon) + 1 = c(Y) +1$ tubes slid into $Y$ to perform our desired isotopy. After dual swapping each tube, it is easy to see that we obtain $S' \cap Y$ as resulting from amalgamation of a stabilized Heegaard splitting of $Y$ (see Figure~\ref{fig:verticalexample}). Hence there is an additional tube which is trivial due to the stabilization, and conclusion $(1)$ of the claim holds. 

\begin{figure}[htbp] 
   \centering
   \includegraphics[width=4in]{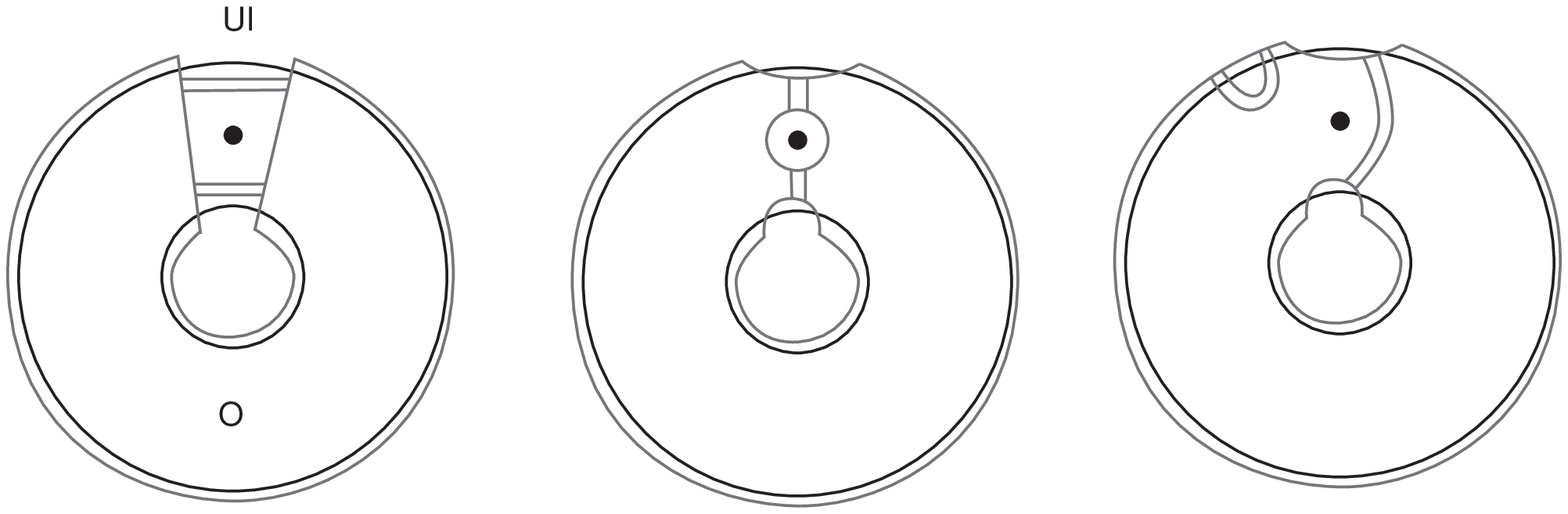} 
   \caption{After dual tube swap we obtain a Heegaard splitting of genus 1 greater than the Heegaard genus, yielding a trivial tube.}
   \label{fig:verticalexample}
\end{figure}

Finally, suppose that $k=0$ and $c(\Upsilon) = -1$ so that $\Upsilon$ contains no exceptional fibers or components of $\partial M$ and none of the $U_i$, $1 \leq i \leq r$, contain exceptional fibers. Either $m_E = 0$, or $m_E \neq 0$. If $m_E \neq 0$, the $U_i$ contain components of $\partial M$, and this implies $(Y, \partial_1 Y, \partial_2 Y)$ is weighted. In either case, $c(Y) = \sum_{i=1}^r c(U_i) + c(\Upsilon) = \sum_{i=1}^r c(U_i) - 1$. Hence, the number of tubes needed to be slid into $Y$ to get spanning disks in $Y$ for each component of $W' \cap (\partial Y - \partial M)$ is $\sum_{i=1}^r c(U_i) = c(Y) + 1$. This is conclusion $(2)$ of the claim. 

\end{proof}

To complete the proof of the lemma, it suffices to show that 

$$a(Y) = \left\{ \begin{array} {ll}
	c(Y) & \mbox{if $k \geq 1$, or $m_E \geq 1$ and $(Y, \partial_1 Y, \partial_2 Y)$ is not weighted} \\
	c(Y) + 1 & \mbox{if $k=0$, and $m_E = 0$ or $(Y, \partial_1 Y, \partial_2 Y)$ is weighted} 
	\end{array}	\right. $$

Showing this implies that $t = a(Y)$ unless conclusion $(1)$ of the claim holds, in which case $t = a(Y) + 1$ and after sliding in $t$ tubes there is a trivial tube in $Y$. 	
\medskip
	
First assume that $k \geq 1$ or $m_E \geq 1$ and $(Y, \partial_1 Y, \partial_2 Y)$ is not weighted. Then,
$$c(Y) = \sum_{i=1}^r c(U_i) + \sum_{i=1}^s c(\Upsilon_i)$$

$$ = \sum_{i=1}^r (|U_i \cap (\partial Y - \partial M)| + \varepsilon(U_i) - 1) + \sum_{i=1}^s (\varepsilon( \Upsilon_i) - 1).$$

\noindent Noting that $| U_i \cap (\partial Y - \partial M)| $ = $|U_i \cap S|$, this becomes

$$= \sum_{i=1}^r (|U_i \cap S| + \varepsilon(U_i) - 1) + \sum_{i=1}^s (\varepsilon(\Upsilon_i) - 1)$$

$$ = | S \cap Y | + m_E + k - r - s$$

\noindent since $\sum_{i=1}^r \varepsilon(U_i) + \sum_{i=1}^s \varepsilon(\Upsilon_i) = k + m_E$. 

Now, if $B$ is the base of $Y$, then an Euler characteristic calculation shows that $\chi(B) = (r + s) - |S \cap Y| - m_E$. The above then becomes

$$ = - \chi(B) + k$$

$$ = -(2 - 2g - m) + k $$

$$ = 2g + m + k - 2$$

Since $ k \geq 1$ or $m_E \geq 1$ and $(Y, \partial_1 Y, \partial_2 Y)$ is not weighted, by Remark~\ref{amalgamationgenusofseifertfiberedspaces} $g_h (Y) = 2g + m + k - 1$. Hence $$2g + m + k - 2 = g_h(Y) - 1 = a(Y).$$

If $k=0$ and either $m_E = 0$ or $(Y, \partial_1 Y, \partial_2 Y)$ is weighted, then as in the proof of the claim, we can assume that $s=1$ and $c(\Upsilon) = - 1$. Thus we have

$$c(Y) = \sum_{i=1}^r c(U_i) + c(\Upsilon)$$

$$= \sum_{i=1}^r (|U_i \cap S| + \varepsilon(U_i) - 1) - 1$$

$$ = |S \cap Y| + m_E - r - 1$$

$$ = - \chi(B) $$

$$= 2g + m - 2.$$

\noindent By Remark~\ref{amalgamationgenusofseifertfiberedspaces}, in this case $g_h (Y) = 2g + m$, and hence it follows that $$c(Y) + 1 = g_h (Y) -1 = a(Y).$$
\end{proof}

Note that the proof of Lemma~\ref{verticalimportlemma} implies that a tube isotoped into $Y$ either is trivial or contributes to the $a(Y)$ tubes needed for each component of $W' \cap (\partial Y - \partial M)$ to be blocked in $Y$. Also, if more than $a(Y)$ tubes are isotoped into $Y$, then there must be a trivial tube which is not needed for every component of $W' \cap (\partial Y - \partial M)$ to have a spanning disk in $Y$. This trivial tube can thus be isotoped freely in $M$ by Lemma~\ref{trivialtubelemma} without changing the conclusion of Lemma~\ref{verticalimportlemma}. 

\section{Proof of the main results}
\label{sec:proofofthemainresults}

In this section we prove Theorem~\ref{maintheorem1} and the generalization of Corollary~\ref{maincorollary} that it implies. 

\subsection{Isotoping tubes to other horizontal components}

Let $X_1, \ldots, X_h$ be the horizontal components of $M$ cut along $\Theta'$ (include the active component $N$ in this list if $N$ is a Seifert fibered component and $S \cap N$ is psuedohorizontal), $Y_1, \ldots, Y_v$ be the vertical components, and $Q_1, \ldots, Q_q$ be the $T^2 \times I$ components. 

\begin{definition}
\label{idealcomponentsdefinition}
\rm Assume that $X_1, \ldots, X_j$ are contained in $M_1$ and $X_{j+1}, \ldots, X_h$ are contained in $M_2$, where $M_1$ and $M_2$ are the manifolds given by $\Theta'$ as in Definition~\ref{mutuallyseparatingdefinition}. Then the components $Z_1, \ldots, Z_n$ of $(\cup_{i=1}^j X_i \cap V) \cup (\cup_{i={j+1}}^h X_i \cap W)$ are the components of $\cup_{i=1}^h X_i$ cut along $S$ requiring tubes to be added to them, as discussed in Remark~\ref{startingstabilizationremark}. We will call $Z_1, \ldots, Z_n$ the {\em ideal components} of the horizontal components, as they correspond to the ideal partition of $(M, \partial_1 M, \partial_2 M)$.
\end{definition}

\begin{lemma}
\label{tubestohorizontalcomponentslemma}
Let $Z_1, \ldots, Z_n$ be the ideal components of $W \cap \cup_{i=1}^h X_i$ and $V \cap \cup_{i=1}^h X_i$. Let $g$ be the genus of \hs, and suppose that $g \geq a(M, \partial_1 M, \partial_2 M)$. Then after one stabilization of \hs, for every $1 \leq i \leq n$, $Z_i$ has some tube $\tau_i$ isotopic into it such that $\tau_i$ is adjacent to a component of $Z_i \cap \Theta'$. 
\end{lemma}

\begin{proof}

Fix some $1 \leq i \leq n$. Suppose that the component $Z_i$ is contained in $V \cap X_{l_i}$. By adding a tube to $W$ in $Z_i$ such that it is adjacent to a component of $Z_i \cap \Theta'$ and then isotoping $S' \cap Z_i$ into tube position, we obtain $t_i$ tubes (note that if $Z_i$ and $Z_j$ both lie in the same horizontal component, then $t_i = t_j$). Applying Lemma~\ref{tubemigrationlemma}, after tube sliding and further isotopy, these tubes can be isotoped into various components of $M$ cut along $\Theta'$ cut along $S$. Note that isotoping a tube out of $Z_i$ into some component $U$ of $M$ cut along $\Theta'$ cut along $S$ leaves all the other components of $M$ cut along $\Theta'$ cut along $S$ unaffected by the isotopy. Set 

$$t_i = x_i + y_i + q_i + p_i$$

\noindent where after isotopy, $x_i$ is the number of tubes that remain in $X_{l_i}$, $y_i$ is the number of tubes isotoped into vertical components, $q_i$ is the number of tubes isotoped into components homeomorphic to $T^2 \times I$, and $p_i$ is the number of tubes isotoped into other components $Z_j$, $j \neq i$. Assume that this is done so that $p_i$ is maximal; in particular any trivial tubes will be isotoped into some $Z_j$ via Lemma~\ref{trivialtubelemma}. Moreover, if two tubes are isotoped into the same $Z_j$, then one of the tubes is trivial. It is straightforward to see that any tube slid out of $Z_i$ must be in one of the 4 categories above, and moreover any spanning disks of components of $\Theta'$ cut along $S'$ formed by tubes can be assumed to be in the correct compression body, as per Remark~\ref{startingstabilizationremark}.
\medskip

Now stabilize \hs \ in $Z_1$ so that the stabilization is adjacent to a component of $Z_1 \cap \Theta'$. 
Then via isotoping tubes using spanning disks, applying Lemmas~\ref{trivialtubelemma}, \ref{tubemigrationlemma}, \ref{T2timesIlemma} and \ref{verticalimportlemma}, isotope tubes into as many $Z_j$ components as possible so that the tubes are adjacent to components of $Z_j \cap \Theta'$. Once a tube is isotoped into some $Z_j$, $S' \cap Z_j$ can be isotoped into tube position and then via tube slides and isotopy of tubes along spanning disks, those tubes can be isotoped to other components of $M$ cut along $\Theta'$ cut along $S$. Reorder the components $Z_1, \ldots, Z_n$ so that $Z_1, \ldots, Z_i$ denote all the components which have a tube isotoped into them via the above process, $p_1 \geq \ldots \geq p_i$, and $i$ is maximal over all choices of isotopy. Note that by Lemma~\ref{trivialtubelemma} we can initially stabilize \hs \ in any of the components $Z_1, \ldots, Z_i$ to obtain the above isotopy. 

Suppose that $i < n$, so there is some component $Z_{i+1}$ such that no tube is isotoped into it after applying the above isotopy. Then the components $Z_{i+1}, \ldots, Z_n$ do not have any tubes isotoped into them from this isotopy, and it must be the case that $p_i = 0$. In particular, there must be no trivial tubes left over after the isotopy. 

Now stabilize $S'$ in $Z_{i+1}$ as well, and isotope as before. Assume that after doing this, all of the components $Z_{i+1}, \ldots, Z_n$ have had a tube isotoped into them and, reordering if necessary, $p_{i+1} \geq \ldots \geq p_n$. The case that more subsets of the $Z_j$ components are needed for this process is more cumbersome and does not require any new arguments beyond the case that there are two such subsets. Hence we prove the lemma for the case of two subsets $Z_1, \ldots, Z_i$ and $Z_{i+1}, \ldots, Z_n$, although the same arguments can be applied to the case when there are more than two subsets. 

Note that $p_n = 0$, for otherwise there would be some tube in $Z_n$ which could be isotoped into another component $Z_j$. If $j \geq i+1$, then the tube is trivial. Thus the tube can be isotoped into some component $Z_j$ where $1 \leq j \leq i$. But then stabilizing \hs \ in $Z_{i+1}$ initially would not only isotope tubes into each of $Z_{i+1}, \ldots, Z_n$, but by isotoping the tube from $Z_n$ into, say $Z_1$ and then performing the initial isotopy, the components $Z_1, \ldots, Z_i$ would also have tubes isotoped into them as well.  This contradicts our assumption that $i$ is maximal. 

Now let $U$ be a component of $M$ cut along $\Theta'$ cut along $S$ which is a solid torus and does not contain an exceptional fiber or a component of $\partial M$. Thus $U$ is contained in a vertical component, a $T^2 \times I$ component, or the active component $N$ if $N$ is an aligned edge manifold. We will call such a component a {\em simple component}. 
\bigskip

\noindent{\bf Claim:} After stabilizing \hs \ in both $Z_1$ and $Z_{i+1}$ as above, there is no simple component $U$ such that each component of $\partial U \cap \Theta'$ has a distinct tube adjacent to it in $U$.   
\bigskip

\begin{proof}
Suppose that after stabilizing \hs \ in $Z_1$ and $Z_{i+1}$, $U$ is a simple component such that each component of $\partial U \cap \Theta'$ has a distinct tube adjacent to it in $U$. 
First note that if all the tubes in $U$ are isotoped from stabilizing \hs \ in $Z_1$, say, then automatically one of the tubes is trivial. This tube can then be isotoped into $Z_{i+1}$ via Lemma~\ref{trivialtubelemma}, a contradiction. Let $\tau_1$ be a tube resulting from stabilizing \hs \ in $Z_1$, and $\tau_2$ be a tube resulting from stabilizing \hs \ in $Z_{i+1}$ such that $\tau_1$ and $\tau_2$ are isotoped into $U$ and each have a foot on some component of $S \cap U$. Let $A_1$ and $A_2$ be the components of $\partial U \cap \Theta'$ to which $\tau_1$ and $\tau_2$ are respectively adjacent. 

\begin{figure}[h] 
   \centering
   \includegraphics[width= 1.5in]{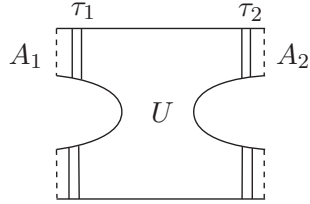} 
   \caption{Each component of $U$ has a distinct tube adjacent to it.}
   \label{fig:setupofu}
\end{figure}

Suppose that $\tau_2$ was isotoped via a tube slide and isotopy along spanning disks to be adjacent to $A_2$. By reversing the the tube slide, we see that $A_2$ has a spanning disk that runs along $\tau_2$ once (see Figure~\ref{fig:runtaurun}). 
\begin{figure}[h] 
   \centering
   \includegraphics[width= 2 in]{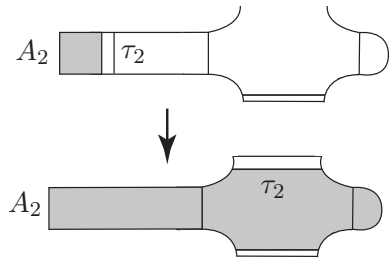} 
   \caption{A spanning disk for $A_2$ running along $\tau_2$ once.}
   \label{fig:runtaurun}
\end{figure}
If $\tau_2$ is trivial, then by definition there is a disk meeting the cocore of $\tau_2$ in a single point. That disk can be taken as a spanning disk of $A_2$ after isotoping $\tau_2$ back to its original position. Note that this spanning disk of $A_2$ can be assumed not to intersect $\tau_2$. Finally, if $\tau_2$ is obtained via a dual tube swap, then reversing the dual tube yields a spanning disk for $A_2$ as in Figure~\ref{fig:dualspan}.

\begin{figure}[h] 
   \centering
   \includegraphics[width= 2.75in]{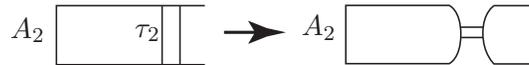} 
   \caption{A spanning disk for $A_2$ obtained by the reversal of a dual tube swap.}
   \label{fig:dualspan}
\end{figure}

The upshot is that for a tube such as $\tau_2$ which, after isotopy, forms a spanning disk for some component of $\Theta'$ cut along $S$ such as $A_2$, then there is a spanning disk for $A_2$ which runs along $\tau_2$ at most one time {\em before} tube sliding and isotopy. 

Let $D_2$ be a spanning disk of $A_2$ which runs along the pre-isotoped $\tau_2$ at most once. In particular $D_2$ can be taken as a spanning disk for $A_2$ before stabilizing \hs \ in $Z_{i+1}$. Thus assume that \hs \ is only stabilized once in $Z_1$. If the disk $D_2$ is disjoint from $\tau_1$, $D_2$ guides an isotopy of $\tau_1$ across $A_2$, as indicated in Figure~\ref{fig:firstisotopyoftau}.

\begin{figure}[h] 
   \centering
   \includegraphics[width= 3.5in]{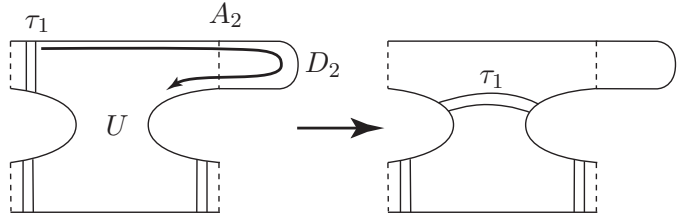} 
   \caption{Isotoping $\tau_1$ across $A_1$ via the disk $D_2$.}
   \label{fig:firstisotopyoftau}
\end{figure}

Note that in Figure~\ref{fig:firstisotopyoftau}, $D_2$ was assumed to be disjoint from the interior of $U$. This can always be done, for otherwise $\tau_2$ would be isotoped into $U$ by going through some other component of $\partial U \cap \Theta'$ first, and we would have taken $\tau_2$ as adjacent to that component. (If $\tau_1$ and $\tau_2$ are adjacent to the same component of $U \cap \Theta'$, then the same arguments apply, taking $U$ as a solid torus between $\tau_1$ and $\tau_2$, and using the disk $D_1$ to isotope $\tau_2$ to be trivial if $D_2$ is not disjoint from the interior of $U$). 

By a symmetric argument, if a spanning disk $D_1$ of $A_1$ before stabilizing \hs \ in $Z_1$ is disjoint from $\tau_2$, then $\tau_2$ can be isotoped across $A_1$.

Now suppose that $D_1$ and $D_2$ are not  disjoint from $\tau_2$ and $\tau_1$, respectively. This implies $D_1$ runs along $\tau_2$ at least once and $D_2$ runs along $\tau_1$ at least once. As this can only occur if both $\tau_1$ and $\tau_2$ are isotoped via tube slides (the spanning disks for $A_1$ and $A_2$ formed by trivial tube isotopies and dual tube swaps can be taken as disjoint from $\tau_1$ or $\tau_2$), it must be the case that $D_1$ and $D_2$ run along $\tau_2$ and $\tau_1$ before stabilization and isotopy via tube slides. This cannot happen, however, as stabilizing \hs \ in $Z_1$ and isotoping $\tau_1$ to be adjacent to $A_1$ would not allow then for $\tau_2$ to be isotoped to be adjacent to $A_2$ after stabilizing \hhs \ additionally in $Z_{i+1}$, as $\tau_2$ then has no disk to be used for tube sliding (see Figure~\ref{fig:secondisotopyoftau}). 

\begin{figure}[h] 
   \centering
   \includegraphics[width= 3.25in]{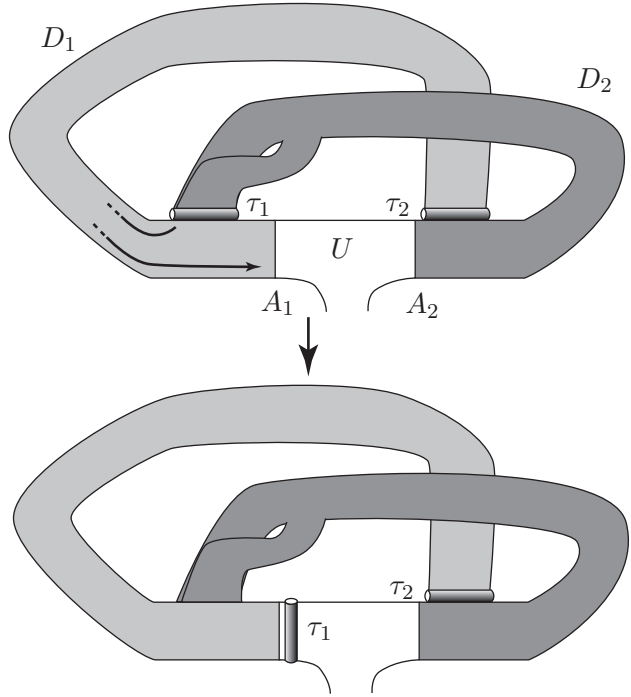} 
   \caption{If $D_1$ and $D_2$ are not disjoint from $\tau_2$ and $\tau_1$ respectively before isotopy, then $\tau_2$ cannot be isotoped to be adjacent to $A_2$ after stabilizing in both $Z_1$ and $Z_{i+1}$.}
   \label{fig:secondisotopyoftau}
\end{figure}

Thus it must be the case that one of $\tau_1$ or $\tau_2$ can be isotoped through $A_2$ or $A_1$, respectively. Repeating this argument for each component of $U \cap \Theta'$, we conclude that one of the tubes isotoped to be adjacent to a component of $U \cap \Theta'$ must be trivial. Without loss of generality, if that tube is one that resulted from stabilizing \hs \ in $Z_1$, then the tube can be isotoped via Lemma~\ref{trivialtubelemma} into $Z_{i+1}$, which is a contradiction. 
\end{proof}

Now, having stabilized \hs \ in both $Z_1$ and $Z_{i+1}$ and isotoping, we have assumed that all components $Z_j$, $1 \leq j \leq n$, have a tube isotopic into them. Moreover, we have $p_{i+1} = p_n = 0$. By the claim, after this isotopy no simple component $U$ is such that every component of $\partial U \cap \Theta'$ contains a distinct tube adjacent to it. 

Suppose that $U$ is in a vertical component. If we can perform a dual tube swap on some tube in $U$ to obtain a trivial tube as in the proof of Lemma~\ref{verticalimportlemma}, then that tube could be isotoped into $Z_1$ or $Z_{i+1}$ and by Lemma~\ref{trivialtubelemma} and we obtain a contradiction. Thus, each simple component contains at most $|\partial U \cap \partial Y| - 1$ tubes, and each non-simple component contains at most $|\partial U \cap \partial Y|$ tubes. Since there are no dual tube swaps yielding trivial tubes, as in the proof of Lemma~\ref{verticalimportlemma} it must be the case that there are at most $\sum_{i=1}^v a(Y_i)$ tubes isotoped into the vertical components, {\em i.e.}, 

$$\sum_{i=1}^n y_i \leq \sum_{i=1}^v a(Y_i).$$

For the components $Q_1, \ldots, Q_q$, Lemma~\ref{T2timesIlemma} implies that at most one tube is isotoped into each $Q_j$. Thus, the number of possible tubes in the simple components of $Q_1, \ldots, Q_q$ is 

$$\sum_{i=1}^n q_i \leq q.$$

Moreover, the if the active component $N$ is an edge manifold, Lemma~\ref{activecomponentlemma} implies that we can assume no tubes are in $N$. If $N$ is Seifert fibered and $S \cap N$ is pseudohorizontal, then we treat it as a horizontal component.

Finally, note that each horizontal component $X_i$ has at most $a(X_i)$ tubes remaining in it after the isotopy of $S'$. Also, since $p_{i-1} = p_n =0$, it must be the case that 

$$\sum_{i=1}^n p_i = n -2.$$

Putting this all together, we have 

$$\sum_{i=1}^n t_i = \sum_{i=1}^n x_i + y_i + q_i +p_i \leq \sum_{i=1}^h a(X_i) + \sum_{i=1}^v a(Y_i) + q + n -2.$$

Now, using an Euler characteristic argument, we have

$$\chi(S) = \sum_{i=1}^h \chi(S \cap X_i) + \sum_{i=1}^v \chi(S \cap Y_i) + \chi(S \cap N).$$

\noindent For each $1 \leq i \leq h$, $S \cap X_i$ is a disjoint union of $2 n_i$ copies of a punctured surface of genus $g_i$ with $m_i$ boundary components. In particular, 

$$\sum_{i=1}^h \chi(S \cap X_i) = \sum_{i=1}^h 2n_i(2-2g_i - m_i)$$

$$ = \sum_{i=1}^h 2n_i - \sum_{i=1}^h 2n_i (2g_i + m_i -1).$$

Now, we note that $\sum_{i=1}^h n_i = n$, the number of $Z_j$ components, and that $\sum_{i=1}^h n_i (2g_i + m_i -1) = \sum_{i=1}^n t_i$, the number of tubes generated by adding a tube to each $Z_j$ component. Using the fact that $S \cap Y_i$ is a disjoint union of annuli for $1 \leq i \leq v$, it must be the case that $\chi(S \cap Y_i) = 0$. The above equation then becomes

$$\chi(S) = 2n - \sum_{i=1}^n 2t_i + \chi(S \cap N).$$

\noindent As \hs \ has genus $g$,

$$\chi(S) = 2 - 2g =  2n - \sum_{i=1}^n 2 t_i + \chi(S \cap N),$$

\noindent which implies

$$2g - 2 = \sum_{i=1}^n 2t_i - 2n - \chi(S \cap N)$$
\noindent and so
$$g = \sum_{i=1}^n t_i  - n - \frac{\chi(S \cap N)}{2} + 1.$$

Above, we established that $\sum_{i=1}^n t_i \leq \sum_{i=1} x_i + y_i + q_i +p_i = \sum_{i=1}^h a(X_i) + \sum_{i=1}^v a(Y_i) + q + n -2.$ If the active component $N$ is an edge manifold, then as $S \cap N$ is a disjoint union of essential annuli and one 4-times punctured sphere component, $-\frac{\chi(S \cap N)}{2} = 1$. Thus 

$$g = \sum_{i=1}^n t_i  - n - \frac{\chi(S \cap N)}{2} + 1 \leq \sum_{i=1}^h a(X_i) + \sum_{i=1}^v a(Y_i) + q + n -2 - n +2$$

$$ = \sum_{i=1}^h a(X_i) + \sum_{i=1}^v a(Y_i) + q  $$

$$ = a(M, \partial_1 M, \partial_2 M) - 1$$

\noindent by Lemma~\ref{graphamalgamationgenus}. As this implies $g < a(M, \partial_1 M, \partial_2 M)$, this is a contradiction.

If $N$ is a Seifert fibered component and $S \cap N$ is pseudohorizontal, then $S \cap N$ is obtained as the connect sum of two copies of a punctured surface of genus $g_N$ with $m_N$ boundary components. Hence $$- \frac{\chi(S \cap N)}{2} = - \frac{2(2-2g_N - m_N) - 2}{2} = -(1- 2g_N - m_N) = t_N.$$ Moreover, in this case $N$ is considered as a horizontal component and the computation of its amalgamation genus can be taken as one of the $a(X_i)$. Thus it follows that 

$$g = \sum_{i=1}^n t_i  - n - \frac{\chi(S \cap N)}{2} + 1 \leq  \sum_{i=1}^h a(X_i) + \sum_{i=1}^v a(Y_i) + q + n -2 - n + 1$$

$$= a(M, \partial_1 M, \partial_2 M) - 2$$

\noindent which, as before, is a contradiction. 

Thus it must be the case that every $Z_j$ component has a tube isotopic into it after at most one stabilization of \hs.
\medskip

The same argument applies in the case that stabilizations are needed for more than two components $Z_1$ and $Z_{i+1}$. Indeed, if after stabilizing \hs \ in $Z_1$ and $Z_{i+1}$ there are still $Z_j$ components which do not have a tube isotoped into them, then we form another set of $Z_j$ components as needed. The same argument as above then yields the result.

\end{proof}

\subsection{Proof of Theorem~\ref{maintheorem1} and Corollary~\ref{maincorollary1}}
The proof of Theorem~\ref{maintheorem1} follows readily.
\bigskip

\noindent{\em Proof of Theorem~\ref{maintheorem1}.}
Suppose $g \geq a(M, \partial_1 M, \partial_2 M)$. As $\Theta'$ is mutually separating, $M$ cut along $\Theta'$ consists of two manifolds $M_1$ and $M_2$ such that a regular neighborhood of each component of $\Theta'$ meets both $M_1$ and $M_2$ nontrivially. By Lemma~\ref{nohorizontalcomponentslemma}, we can assume that some component of $M$ cut along $\Theta'$ is horizontal, or that $N$ is a Seifert fibered component and $S \cap N$ is pseudohorizontal. Thus, by Lemma~\ref{tubestohorizontalcomponentslemma}, stabilizing \hs \ once yields a tube isotopic into each of the ideal components $Z_1, \ldots, Z_n$. This yields $$\sum_{i=1}^n t_i  - \frac{\chi(S \cap N)}{2} = g + (n-1) $$ 
$$\geq a(M, \partial_1 M, \partial_2 M) + (n-1)$$
$$= \sum_{i=1}^h a(X_i) + \sum_{i=1}^v a(Y_i) +1  + q + (n-1) $$tubes by Lemma~\ref{graphamalgamationgenus}. By Lemmas~\ref{tubemigrationlemma}, \ref{T2timesIlemma}, \ref{activecomponentlemma} and \ref{verticalimportlemma}, tubes can be isotoped into each vertical and $T^2 \times I$ component so that every component of $V' \cap \Theta'$, say, has a spanning disk in $M_1$ and every component of $W' \cap \Theta'$ has a spanning disk in $M_2$. By Lemma~\ref{theamalgamationlemma}, \hhs \ is an amalgamation along $\Theta'$.

As \hhs \ is a stabilized amalgamation along $\Theta'$, by Corollary~\ref{graphmanifoldsamalgamationlemma}, \hhs \ and any other amalgamation along $\Theta'$ are isotopic after at most one stabilization.
\medskip

Now assume $g < a(M, \partial_1 M, \partial_2 M)$. In this situation, stabilize \hs \ $a(M, \partial_1 M, \partial_2 M) - g + 1$ times so that the genus of \hhs \ is $a(M, \partial_1 M, \partial_2 M) + 1$. The above arguments can then be applied to show that \hhs \ is an amalgamation along $\Theta'$, and hence isotopic to any other such amalgamation along $\Theta'$ after at most one stabilization of the other splitting.

\begin{flushright} 
$\Box$
\end{flushright}

\begin{corollary}
\label{maincorollary1}
If $(M, \partial_1 M, \partial_2 M)$ has two strongly irreducible Heegaard splittings, at least one of which has genus greater than or equal to $a(M, \partial_1 M, \partial_2 M)$, then the two splittings are isotopic after at most one stabilization of the larger genus splitting. 
\end{corollary}

\begin{proof}
Let \hs \ and \pq \ be strongly irreducible Heegaard splittings of $(M, \partial_1 M, \partial_2 M)$ such that one of $S$ or $\Sigma$ has genus greater than or equal to $a(M, \partial_1 M, \partial_2 M)$. By stabilizing each splitting the appropriate number of times (only once for the higher genus splitting), \hs \ stabilizes to \hhs \ and \pq \ stabilizes to \ppq . By Theorem~\ref{maintheorem1}, both \hhs \ and \ppq \ are (stabilized) amalgamations along $\Theta'$. By Corollary~\ref{graphmanifoldsamalgamationlemma}, \hhs \ and \ppq \ are isotopic. 
\end{proof}
 
\section{Examples of small and large genus Heegaard splittings}
\label{sec:section11}

Theorem~\ref{maintheorem1} implies that the genus of a strongly irreducible Heegaard splitting \hs \ needs to at least as big as $a(M, \partial_1 M, \partial_2 M)$ to be isotopic to an amalgamation along $\Theta'$. This raises the natural question, When is the genus of \hs \ at least as big as $a(M, \partial_1 M, \partial_2 M)$? In what follows, we provide examples of closed graph manifolds $M$ having genus $g$ Heegaard splittings \hs \ of $M$ with $g < a(M)$ and with $g \geq a(M)$. In the latter case, we show also that the constructed splittings are strongly irreducible. 

\subsection{Small genus splittings} 

Let $X$ be a Seifert fibered space with Seifert data 

$$\left( 0; m; \frac{1}{3}, \frac{2}{3} \right).$$ 

Then there exists a connected horizontal surface $S_0 \cap X$ in $X$ such that $|\partial S_0| = 3m$. This surface is a planar surface with $3m$ boundary components. Take $S \cap X$ as two parallel copies of $S_0$ in $X$. Thus $X$ cut along $S$ is a genus $3m - 1$ handlebody. For $m-1$ of the boundary components of $X$, attach a copy of a Seifert fibered space $Y$ with Seifert data 

$$\left( 0; 1; \frac{\beta_1}{\alpha_1}, \frac{\beta_2}{\alpha_2}, \frac{\beta_3}{\alpha_3}, \frac{\beta_4}{\alpha_4} \right)$$

\noindent such that $S \cap Y$ is a disjoint union of 3 vertical annuli separating the exceptional fibers of $Y$, and $S \cap Y$ is glued to $S \cap X$. 

Let $N$ be homeomorphic to $T^2 \times I$ and take a disjoint union of five essential annuli, exactly four of which have their boundary components on different components of $\partial N$. Let $S \cap N$ be the surface obtained from the four essential annuli in $N$ via performing ambient 1-surgery along an arc which connects the annulus with both boundary components on the same component of $\partial N$ to one of the other annuli (so $S \cap N$ is aligned). Let $T$ be the component of $\partial N$ such that $| S \cap T| = 6$. 

Glue the remaining component of $\partial X$ to $T$ such that $S \cap X$ is glued to $S \cap N$. To the other component of $\partial N$, glue a Seifert fibered space $Y'$ with Seifert data 

$$\left( 0; 1; \frac{\beta_1}{\alpha_1}, \frac{\beta_2}{\alpha_2} , \frac{\beta_3}{\alpha_3} \right)$$ 

\noindent such that $S \cap Y'$ is two vertical annuli separating the exceptional fibers of $Y'$, and $S \cap Y'$ is glued to $S \cap N$. Thus the manifold obtained via gluing is a closed graph manifold $M$ with $m +1$ Seifert fibered components. It is a nice exercise to show that $M$ cut along $S$ is two handlebodies $V$ and $W$, so that \hs \ is a Heegaard splitting of $M$. Moreover, it is a simple calculation to show that the genus of \hs \ is $3m$. By Lemma~\ref{graphamalgamationgenus}, $a(M) = a(X) + m(a(Y)) + a(Y') + 1 = m + 3m + 2 + 1 = 4m + 3$. Hence $g = 3m < 4m+3 = a(M)$. 

These examples can be generalized by modifying the numbers of exceptional fibers in $X$, $Y$ and $Y'$. Note that what makes $g < a(M)$ is that fact that $S \cap X$ is has low genus. In general, for $S_0$ a component of $S \cap X$, $\chi(S_0) = 2 - 2g_0 - m_0$ where $g_0$ is the genus of $S_0$ and $m_0$ is the number of boundary components. Note that $\chi(S_0) = 1 - t$, where $t$ is the number of tubes generated by isotoping $S' \cap X$ to tube position after stabilization. 

In the above examples, $g_0 = 0$, and each component of $\partial S_0$ required a tube slide to become an amalgamation along $\Theta$. Since $2 - m_0 = 1 - t$, we have $t = m_0 - 1$, and since $a(X)$ of the tubes must remain in $X$, there are only $t - a(X)$ available for tube sliding, a quantity clearly less than $m_0$. In particular, if $2g_0 < a(X) + 1$, then $t - a(X) < m_0$ which implies $g < a(M)$ in the above examples. 

\subsection{Large genus splittings}

The following examples are due to D. Bachman. Let $X_1$ be a torus knot complement, hence a Seifert fibered space over a disk with 2 exceptional fibers of multiplicity $p$ and $q$. Let $S \cap X_1$ be a bridge surface in $X_1$, hence $S \cap X_2$ is a $(p-1)(q-1)$-times punctured sphere. Let $X_2$ be any Seifert fibered space with Seifert data $\left(g;1;\frac{\beta_1}{\alpha_1}, \ldots, \frac{\beta_k}{\alpha_k}\right)$ containing a horizontal surface $F$ with one boundary component. Without loss of generality, assume $k \geq 1$. Let $S \cap X_2$ be $(p-1)(q-1)$ copies of $F$. Glue $X_1$ to $X_2$ so that $S \cap X_1$ is glued to $S \cap X_2$. This gives a graph manifold $M$. It is a nice exercise to show that $S$ gives a Heegaard splitting \hs \ of $M$. 

The Heegaard splitting \hs \ is strongly irreducible since for any torus knot, bridge position equals thin position (see {\em e.g.}~\cite{Thompson}). That is, since $X_2$ cut along $S$ consists of $(p-1)(q-1)$ copies of $F \times I$, any compressing disks of $V$ or $W$ must have outermost disk components in $X_1$. If two outermost disk components $D$ and $E$ in $V \cap X_1$ and $W \cap X_1$ respectively are disjoint, then this would define an isotopy of the torus knot in $X_1$ showing that thin position is not bridge position, a contradiction. Hence every pair of disks, one in $V$ and one in $W$, must interesect in $X_1$, implying \hs \ is strongly irreducible. 

One can check that $a(M) =  2g + k$. Moreover, as $F$ branch covers the base of $X_2$, the multiplicity of the cover is divisible by $\mbox{lcm} (\alpha_1, \ldots, \alpha_k)$. By choosing multiplicities of sufficient size, the Euler characteristic of $F$, and hence the genus of $S$, may be chosen to be arbitrarily high. This implies $g > a(M)$. Moreover, if $X_2$ contains any vertical incompressible tori (a condition satisfied, for example, if $k \geq 3$ or $g \geq 1$) then Dehn twisting along these tori produces other strongly irreducible Heegaard splittings of the same genus of $M$. Corollary~\ref{maincorollary} implies that these splittings are isotopic after at most one stabilization.

\bigskip

\bigskip

\noindent \address{Department of Mathematics \\ The University of Texas at Austin \\ Austin, TX 78712 USA}
\bigskip

\noindent \email{rdtalbot@math.utexas.edu}

\end{document}